# BACKWARD STOCHASTIC DIFFERENTIAL EQUATIONS WITH RANDOM STOPPING TIME AND SINGULAR FINAL CONDITION

By A. Popier

*Université de Provence*

In this paper we are concerned with one-dimensional backward stochastic differential equations (BSDE in short) of the following type:

$$Y_t = \xi - \int_{t\wedge\tau}^{\tau} Y_r |Y_r|^q \, dr - \int_{t\wedge\tau}^{\tau} Z_r \, dB_r, \qquad t \geq 0,$$

where $\tau$ is a stopping time, $q$ is a positive constant and $\xi$ is a $\mathcal{F}_\tau$-measurable random variable such that $\mathbf{P}(\xi = +\infty) > 0$. We study the link between these BSDE and the Dirichlet problem on a domain $D \subset \mathbb{R}^d$ and with boundary condition $g$, with $g = +\infty$ on a set of positive Lebesgue measure.

We also extend our results for more general BSDE.

**Introduction.** Let $(\Omega, \mathcal{F}, \mathbf{P})$ be a probability space, $B = (B_t)_{t\geq 0}$ a Brownian motion defined on this space, with values in $\mathbb{R}^d$. $(\mathcal{F}_t)_{t\geq 0}$ is the standard filtration of the Brownian motion. Also given are $\tau$ a $\{\mathcal{F}_t\}$-stopping time, $\xi$ a real, $\mathcal{F}_\tau$-measurable random variable, called the *final condition*, and $f : \Omega \times \mathbb{R}^+ \times \mathbb{R} \times \mathbb{R}^d \to \mathbb{R}$ the *generator*.

We wish to find a progressively measurable solution $(Y, Z)$, with values in $\mathbb{R} \times \mathbb{R}^d$, of the BSDE

(1) $$Y_t = \xi + \int_{t\wedge\tau}^{\tau} f(r, Y_r, Z_r) \, dr - \int_{t\wedge\tau}^{\tau} Z_r \, dB_r, \qquad t \geq 0.$$

Such equations, in the nonlinear case, have been introduced by Pardoux and Peng in 1990 in [19], when $\tau$ is replaced by a constant time $T > 0$. They gave the first existence and uniqueness result. Since then, BSDE have been studied with great interest (see the references in [18]). In particular, Peng [20] describes how the solution $Y$ of (1) for an unbounded random terminal time is related to a semilinear elliptic PDE. Viscosity solutions for such equations will be constructed by stochastic methods (see Theorem 8 below). This generalization of the Feynman–Kac formula is a reason for studying random terminal times.

Let us recall the definition of a solution of (1) which can be found in [4].

DEFINITION 1. A solution of the BSDE (1) is a pair $\{(Y_t, Z_t), t \geq 0\}$ of progressively measurable processes with values in $\mathbb{R} \times \mathbb{R}^d$ such that, **P**-a.s.:

- on the set $\{t \geq \tau\}$, $Y_t = \xi$ and $Z_t = 0$,
- $t \mapsto \mathbf{1}_{t \leq \tau} f(t, Y_t, Z_t)$ belongs to $L^1_{\text{loc}}(0, \infty)$, $t \mapsto Z_t$ belongs to $L^2_{\text{loc}}(0, \infty)$,
- and for all $0 \leq t \leq T$,

$$Y_{t \wedge \tau} = Y_{T \wedge \tau} + \int_{t \wedge \tau}^{T \wedge \tau} f(r, Y_r, Z_r)\, dr - \int_{t \wedge \tau}^{T \wedge \tau} Z_r\, dB_r.$$

A solution is said to be an $L^p$-solution for some $p > 1$ if, moreover, for some $\lambda \in \mathbb{R}$,

$$\mathbb{E}\left(\sup_{0 \leq t \leq \tau} e^{p\lambda t}|Y_t|^p + \int_0^\tau e^{p\lambda t}|Y_t|^p\, dt + \int_0^\tau e^{p\lambda t}|Y_t|^{p-2}\|Z_t\|^2\, dt\right) < +\infty.$$

We assume that the generator $f : \Omega \times \mathbb{R}^+ \times \mathbb{R} \times \mathbb{R}^d \to \mathbb{R}$ is such that:

(H0) $f(\cdot, y, z)$ is progressively measurable, for all $y, z$;

(H1) $\exists K \geq 0$, such that a.s. $\forall t, y, z, z'$,

$$|f(t, y, z) - f(t, y, z')| \leq K\|z - z'\|;$$

(H2) $\exists \mu \in \mathbb{R}$, such that a.s. $\forall t, y, y', z$,

$$(y - y')(f(t, y, z) - f(t, y', z)) \leq \mu|y - y'|^2;$$

(H3) $y \mapsto f(t, y, z)$ is continuous, $\forall t, z$, a.s.

(H4) for all $r > 0$ and all $n \in \mathbb{N}^*$, $\psi_r(t) = \sup_{|y| \leq r} |f(t, y, 0) - f(t, 0, 0)|$ belongs to $L^1((0, n) \times \Omega)$.

Now for some $p > 1$ we suppose that there exists $\lambda > \nu_p = \mu + \frac{K^2}{2(p-1)}$, such that

(H5) $$\mathbb{E}\left[\int_0^\tau e^{p\lambda t}|f(t, 0, 0)|^p\, dt\right] < +\infty$$

and

(H6) $$\mathbb{E}\left[e^{p\lambda \tau}|\xi|^p + \int_0^\tau e^{p\lambda t}|f(t, e^{-\nu_p t}\overline{\xi}_t, e^{-\nu_p t}\overline{\eta}_t)|^p\, dt\right] < +\infty,$$

where $\overline{\xi} = e^{\nu_p \tau}\xi$, $\overline{\xi}_t = \mathbb{E}(\overline{\xi}|\mathcal{F}_t)$ and $\overline{\eta}$ is predictable and such that

$$\overline{\xi} = \mathbb{E}(\overline{\xi}) + \int_0^{+\infty} \overline{\eta}_t\, dB_t, \qquad \mathbb{E}\left[\left(\int_0^\infty |\overline{\eta}_t|^2\, dt\right)^{p/2}\right] < \infty.$$

Let us recall Theorem 5.2 of [4].

THEOREM 1. *Under the conditions* (H0)–(H6), *there exists a unique solution* $(Y,Z)$ *of the BSDE* (1), *which, moreover, satisfies, for* $\lambda > \nu_p$ *such that* (H5) *and* (H6) *hold:*

$$
\begin{aligned}
\mathbb{E}&\bigg(\sup_{0\leq t\leq \tau} e^{p\lambda t}|Y_t|^p + \int_0^\tau e^{p\lambda r}|Y_r|^{p-2}(|Y_r|^2 + \|Z_r\|^2)\,dr\bigg)\\
&\qquad \leq c\mathbb{E}\bigg(e^{p\lambda \tau}|\xi|^p + \int_0^\tau e^{p\lambda r}|f(t,0,0)|^p\,dr\bigg),
\end{aligned}
\tag{2}
$$

*for some constant* $c = c(p,\lambda,K,\mu)$.

REMARK 1. The previous theorem is a generalization of the result of Darling and Pardoux (Theorem 3.4 in [6]) or of Pardoux (Theorem 4.1 in [18]). In [6] or [18] the result is given in the case $p=2$. Here we have expressed the theorem for the dimension one ($\xi$ and $Y_t$ belong to $\mathbb{R}$). But it is still true in higher dimensions (see [4]; the product in (H2) must be replaced by the scalar product in $\mathbb{R}^m$).

Note that if $f$ is a Lipschitz function, the condition (H2) holds.

From now and in the rest of the paper we are concerned with the BSDE

$$
Y_t = \xi - \int_{t\wedge\tau}^\tau Y_r|Y_r|^q\,dr - \int_{t\wedge\tau}^\tau Z_r\,dB_r \qquad \text{with } q > 0.
\tag{3}
$$

Here the function $f$ is deterministic and equal to

$$f(t,y,z) = -y|y|^q.$$

$f$ satisfies all conditions (H0)–(H4) of Theorem 1, with $K = \mu = 0$ (which implies $\nu_p = 0$ for all $p > 1$). Indeed, $f$ is a nonincreasing function, thereby,

$$-(y-y')(y|y|^q - y'|y'|^q) \leq 0.$$

Since $f(t,0,0) \equiv 0$, (H5) is always satisfied.

The stopping time $\tau$ is defined as follows. Let $D$ be an open bounded subset of $\mathbb{R}^d$, whose boundary is at least of class $C^2$ (see [12] for the definition of a regular boundary). For all $x \in \mathbb{R}^d$, let $X^x$ denote the solution of the SDE:

$$
X_t^x = x + \int_0^t b(X_r^x)\,dr + \int_0^t \sigma(X_r^x)\,dB_r \qquad \text{for } t \geq 0.
\tag{4}
$$

The functions $b$ and $\sigma$ are defined on $\mathbb{R}^d$, with values respectively in $\mathbb{R}^d$ and $\mathbb{R}^{d\times d}$, and are measurable such that:

- Lipschitz condition: there exists $K \geq 0$ such that

$$\forall (x,y) \in \mathbb{R}^d \times \mathbb{R}^d \qquad \|\sigma(x) - \sigma(y)\| \leq K|x-y|; \tag{L}$$

- Boundedness condition:

(B) $$\forall x \in \mathbb{R}^d \qquad |b(x)| + \|\sigma(x)\| \leq K;$$

- Uniform ellipticity: there exists a constant $\alpha > 0$ such that

(E) $$\forall x \in \mathbb{R}^d \qquad \sigma\sigma^*(x) \geq \alpha \mathrm{Id}.$$

In the rest of this paper (L), (B) and (E) are supposed to be satisfied. Under these assumptions, from a result of Yu Veretennikov [24] and [25], equation (4) has a unique strong solution $X^x$. For each $x \in \overline{D}$, we define the stopping time

(5) $$\tau = \tau_x = \inf\{t \geq 0, X_t^x \notin \overline{D}\}.$$

Our stopping time satisfies the following two properties. Since $D$ is bounded and since the conditions (B) and (E) hold

(C1) \qquad\qquad\qquad every point $x \in \partial D$ is regular.

In particular, if $x \in \partial D$, $\tau_x = 0$ a.s. (see [3], Corollary 3.2). This assumption (C1) is important to define a singular solution (see Definition 2 below). Moreover, since (L), (B) and (E) hold, we have the following result (see [21], Theorem 2.1 and [18], Remark 5.6): for all $x \in \overline{D}$, $\tau_x < +\infty$ a.s. and there exists $\beta > 0$ such that

(C2) $$\sup_{x \in \overline{D}} \mathbb{E}(e^{\beta \tau_x}) < \infty.$$

This property will be used to construct solutions of the BSDE (3) for bounded terminal conditions $\xi$ (see Proposition 2).

From the papers [6, 18] or [20], we know that the BSDE (3) with terminal time equal to $\tau = \tau_x$ and final data equal to $\xi = h(X_{\tau_x}^x)$ is associated with the following elliptic PDE with Dirichlet condition $h$:

(6) $$\begin{aligned} -\mathcal{L}u + u|u|^q &= 0 \qquad \text{on } D, \\ u &= h \qquad \text{on } \partial D; \end{aligned}$$

where $\mathcal{L}$ is the second order partial differential operator: for all $\varphi \in C_0^2(\mathbb{R}^d)$,

(7) $$\forall x \in \mathbb{R}^d \qquad \mathcal{L}\varphi(x) = \tfrac{1}{2}\mathrm{Trace}(\sigma\sigma^*(x)D^2\varphi(x)) + b(x)\nabla\varphi(x).$$

In the rest of this paper $\nabla$ and $D^2$ will denote respectively the gradient and the Hessian matrix. If $(Y^x, Z^x)$ denotes the solution of the BSDE (3) with terminal data $h(X_{\tau_x}^x)$, the connection is given by the formula

$$u(x) = Y_0^x.$$

Le Gall [13] succeeded in describing all solutions of the equation $\Delta u = u^2$ in the unit disk $D$ in $\mathbb{R}^2$ by a purely probabilistic method. He established a 1–1 correspondence between all solutions and all pairs $(\Gamma, \nu)$, where $\Gamma$ is a closed subset of $\partial D$ and $\nu$ is a Radon measure on $\partial D \setminus \Gamma$. The set $\Gamma$ is the set of singular points of $\partial D$ where the solution explodes badly: roughly speaking, near points of $\Gamma$, the solution behaves like the inverse of the squared distance to the boundary. The measure $\nu$ can be interpreted as the "boundary value" of $u$ on $\partial D \setminus \Gamma$. The solution corresponding to $(\Gamma, \nu)$ is expressed in terms of the Brownian snake (a path-valued Markov process). In [14] the results announced in [13] are proved in detail and are extended to a general smooth domain in $\mathbb{R}^2$.

The pair $(\Gamma, \nu)$ is called the boundary trace for positive solution of the PDE (6). The definition of boundary trace in general was provided by Marcus and Véron [15] who showed by analytic methods that every positive solution of (6) possesses a unique trace. The trace can be described by a (possibly unbounded) positive regular Borel measure $\tilde{\nu}$ on $\partial D$. The correspondence between $(\Gamma, \nu)$ and $\tilde{\nu}$ is given by

$$\tilde{\nu}(A) = \begin{cases} \nu(A), & \text{if } A \subseteq (\partial D \setminus \Gamma), \\ \infty, & \text{if } A \cap \Gamma \neq \varnothing, \end{cases}$$

for every Borel subset $A$ of $\partial D$.

The corresponding boundary value problem is presented in [15] in the subcritical case $0 < q < 2/(d-1)$ and in [16] in the supercritical case $q \geq 2/(d-1)$. In the subcritical case, for every pair $(\Gamma, \nu)$, the problem has a unique solution. Remark that in [13] and [14], $q = 1$ and $d = 2$, that is, the subcritical case is studied: $q = 1 < 2/(2-1) = 2/(d-1)$. In the supercritical case Marcus and Véron derive necessary and sufficient conditions for the existence of a maximal solution. Similar conditions were obtained by Dynkin and Kuznetsov [8] for $q \leq 1$. Their method relies on probabilistic techniques and is not extendable to $q > 1$, because the main tool is the $q$-superdiffusion which is not defined for $q > 1$.

The object of the present paper is to give a probabilistic representation of the solution of the PDE (6) in terms of the solution of the related BSDE (3). In general, a solution of the PDE has a "blow-up" set $\Gamma$. Therefore, the final data $\xi$ of the BSDE must be allowed to be infinite with positive probability and the set $\{\xi = +\infty\}$ corresponds to $\Gamma$. Hence, our first problem is to find a solution of (3) when $\xi$ is infinite with positive probability, which implies, in particular, that (H6) is not satisfied.

Note that there are some differences between our work and the results of Le Gall or Dynkin and Kuznetsov. With the superprocesses (see [14] or [8]), it should be assumed that $q \leq 1$. In our case there is no restriction on $q > 0$.

Moreover, the Dirichlet boundary condition for the PDE (6) is not taken in the same sense in the two approaches. With the notion of the *boundary trace* (see [8, 14, 15] and [16]), there always exists a maximal positive solution; if $q < 2/(d-1)$, this solution is unique, and if $q \geq 2/(d-1)$, the problem (6) may possess more than one positive solution. More precisely, assume that $D$ is the unit ball in $\mathbb{R}^d$, that $q \geq 2/(d-1)$ and denote by $\mu_\infty$ the Borel measure on $\partial D$ which assigns the value $+\infty$ to every nonempty set. Then for every $\varepsilon > 0$, there exists a positive solution of (6) such that $u(0) < \varepsilon$ and the trace of $u$ is $\mu_\infty$ (see Proposition 5.1 of [16]).

In our case the Dirichlet condition in (6) is taken in the viscosity sense (see Definition 4 in Section 5). The results are rather different: there exists a minimal positive viscosity solution. But we are unable to give conditions to ensure uniqueness of the solution.

*Main results.* In the first section we will prove an a priori estimate which is a probabilistic generalization of the Keller–Osserman inequality.

In Section 2 and in the rest of the paper we assume

$$\xi \geq 0 \quad \text{a.s.}$$

and we allow $\xi$ to be infinite with positive probability: $\mathbf{P}(\xi = +\infty) > 0$. We must modify Definition 1 of a BSDE when $\xi$ does not satisfy the condition (H6).

In the rest of this paper $\rho$ denotes the distance from the boundary of $D$. For $x \in \overline{D}$, for all positive $\eta$, let us define the stopping time

(8) $$\tau_\eta^x = \inf\{t \geq 0, \rho(X_t^x) \leq \eta\}.$$

REMARK 2. For $x \in \overline{D}$, $\tau_\eta^x \leq \tau_x$ a.s. and if $x \in D$, when $\eta$ goes to 0, $\tau_\eta^x$ converges to $\tau_x$ a.s. When $x \in \partial D$, for all $\eta > 0$, $\tau_\eta^x = \tau_x = 0$ a.s., because every point $x \in \partial D$ is regular [condition (C1)].

Therefore, we suppose $x$ to be in $D$ and for convenience, we omit the variable $x$.

DEFINITION 2. For an $\mathcal{F}_\tau$-measurable $\xi$ such that $\mathbf{P}(\xi \geq 0) = 1$ and $\mathbf{P}(\xi = \infty) > 0$, the process $(Y, Z)$ is a solution of the BSDE (3) if:

(D1) for all $\eta > 0$ and all $T \geq 0$,

$$\mathbb{E}\left(\sup_{0 \leq t \leq T} |Y_{t \wedge \tau_\eta}|^2 + \int_0^{T \wedge \tau_\eta} |Z_r|^2 \, dB_r\right) < +\infty;$$

(D2) $\mathbf{P}$-a.s. for all $0 \leq t \leq T$ and all $\eta > 0$,

$$Y_{t \wedge \tau_\eta} = Y_{T \wedge \tau_\eta} - \int_{t \wedge \tau_\eta}^{T \wedge \tau_\eta} Y_r |Y_r|^q \, dr - \int_{t \wedge \tau_\eta}^{T \wedge \tau_\eta} Z_r \, dB_r;$$

(D3) on the set $\{t \geq \tau\}$, $Y_t = \xi$ and $Z_t = 0$, and **P**-a.s.,

$$\lim_{t \to +\infty} Y_{t \wedge \tau} = \xi.$$

We first construct a process $\{(Y_t, Z_t); t \geq 0\}$ satisfying the conditions (D1) and (D2) of the previous definition. $(Y, Z)$ is the limit of the sequence of processes $(Y^n, Z^n)$, solution (in the sense of Definition 1) of the BSDE (3) with terminal condition $\xi \wedge n$. From the first section we already know that there exists a constant $C$ such that

$$\forall t \geq 0 \qquad \rho^{2/q}(X_{t \wedge \tau})Y_t \leq C.$$

Moreover, we prove the following:

PROPOSITION 1. *On $\{\xi = +\infty\}$ the explosion rate of $Y$ is in the order of $\rho^{-2/q}(X_{t \wedge \tau})$: there exists a positive constant $\widetilde{C}$ depending on $D$, $q$, the bound on $b$ and $\sigma$ in* (B) *and on the constant $\alpha$ in* (E)*, such that*

$$\liminf_{t \to +\infty} \rho^{2/q}(X_{t \wedge \tau})Y_{t \wedge \tau} \geq \widetilde{C} \qquad \text{a.s. on } \{\xi = +\infty\}.$$

Without other assumption on $\xi$, we cannot prove that $(Y, Z)$ satisfies the condition (D3) of Definition 2.

In Section 3 we first prove that **P**-a.s. the limit of $Y_{t \wedge \tau}$ as $t$ goes to $+\infty$ exists and

$$\lim_{t \to +\infty} Y_{t \wedge \tau} \geq \xi.$$

Then we add some assumptions on $\xi$ and on the diffusion $X$ to insure that the condition (D3) holds. We prove the following:

THEOREM 2. *Under the assumptions:*

- *the terminal data $\xi$ satisfies*

(A1) $$\xi = g(X_\tau),$$

  *where $g : \mathbb{R}^d \to \overline{\mathbb{R}}_+$ is a function such that $F_\infty = \{g = +\infty\} \cap \partial D$ is a closed set;*
- *on $\mathbb{R}^d \setminus F_\infty$, $g$ is locally bounded, that is, for all compact set $\mathcal{K} \subset \mathbb{R}^d \setminus F_\infty$,*

(A2) $$g \mathbf{1}_\mathcal{K} \in L^\infty(\mathbb{R}^d).$$

- *the boundary $\partial D$ belongs to $C^3$;*

*the process $Y$ is continuous, that is, $\lim_{t \to +\infty} Y_{t \wedge \tau} = \xi$ **P**-a.s.*

In the next section we prove if there exists a solution $(\overline{Y}, \overline{Z})$ of the BSDE (3) in the sense of Definition 2, then $\overline{Y} \geq Y$. Therefore, if the process $(Y, Z)$ is a solution (e.g., if the assumptions of Theorem 2 hold), it is the minimal solution.

In the last section we show the connection between the BSDE (3) with terminal condition $g(X^x_{\tau_x})$ and the PDE (6) with Dirichlet condition $g$. The assumptions of Theorem 2 hold. In the previous sections we have defined a process $\{(Y^x_t, Z^x_t); t \geq\}$ which is a solution of the BSDE (3) with terminal data $g(X^x_{\tau_x})$. Next we define

$$u(x) = Y^x_0.$$

The main result follows:

THEOREM 3. *Under the assumptions of Theorem 2, $u$ is a viscosity solution of the PDE (6) with Dirichlet condition $g$.*

Here we do not suppose that a viscosity solution is continuous. But under some stronger assumptions on the operator $\mathcal{L}$, we also give some regularity properties of the solution $u$. We also prove that $u$ is the minimal solution.

In the last section we will see that these results are still true with more general generators $f$.

THEOREM 4. *Assume that $f$ is a nonincreasing and $C^1$ function with $f(0) = 0$, and such that there exists $q > 0$, $\kappa > 0$ s.t.,*

$$\forall y \geq 0 \qquad f(y) \leq -\kappa y^{1+q}. \tag{9}$$

*If $\xi$ is a nonnegative random variable, with $\mathbf{P}(\xi = +\infty) > 0$, and such that the assumptions of Theorem 2 hold, then there exists a process $(Y, Z)$, solution of the BSDE*

$$Y_t = \xi + \int_{t \wedge \tau}^{\tau} f(Y_r) \, dr - \int_{t \wedge \tau}^{\tau} Z_r \, dB_r \tag{10}$$

*[in the sense of Definition 2, with $f$ instead of $y \mapsto -y|y|^q$ in (D2)].*

Moreover the conclusion of Theorem 3 is still true: there exists a minimal viscosity solution for the PDE

$$\begin{aligned} \mathcal{L}u + f(u) &= 0 \quad \text{on } D, \\ u &= g \quad \text{on } \partial D. \end{aligned} \tag{11}$$

*Important remark on the condition* (E). The condition (E) can be relaxed. In the rest of the paper we can also work with the assumptions (L), (B) and we add the following condition: $b$ is continuous and satisfies the monotonicity condition: there exists $\mu \in \mathbb{R}$ such that

(M) $\quad \forall (x,y) \in \mathbb{R}^d \times \mathbb{R}^d \quad \langle x-y|b(x)-b(y)\rangle \leq \mu |x-y|^2;$

here $\langle \cdot | \cdot \rangle$ denotes the scalar product in $\mathbb{R}^d$. Under these assumptions equation (4) has a unique strong solution $X^x$. For each $x \in \overline{D}$, we define the stopping time

$$\tau = \tau_x = \inf\{t \geq 0, \ X_t^x \notin \overline{D}\}.$$

We also assume that the conditions (C1) and (C2) hold.

Under the assumptions (M), (L), (B), (C1) and (C2), the results which may be false are in Section 2, Proposition 1, and in Section 5, Propositions 11 and 12. We are unable to control the explosion rate of $Y$ (see Remark 5), nor to prove that the viscosity solution $u$ is continuous on $D$ without the ellipticity condition.

In the rest of the paper all results (except maybe Propositions 1, 11 and 12) could be proved without the condition (E). Indeed, we use this assumption only in the proofs of Propositions 4 and 8, in order to control the Green function $G(x, \cdot)$ associated to the process $X^x$ killed at $\tau_x$. Under (E), this function $G(x, \cdot)$ is continuous on $D$ except at the point $x$, and is integrable on $D$. This assumption on $G$ can replace (E) (see, e.g., [21] for more details on $G$).

**1. An a priori estimate.** Let $(Y, Z)$ be the solution of the BSDE (3) with terminal data $\xi$ such that the hypothesis (H6) holds. We will need an a priori inequality in order to control $Y_{t \wedge \tau}$ for $t \in [0, +\infty[$. The idea comes from the Keller–Osserman inequality which is true for any open set $D$ ([11] and [17]). Denote by $\rho$ the distance to the boundary of $D \subset \mathbb{R}^d$.

THEOREM 5 (Keller–Osserman). *There exists a positive constant* $C = C(q,d)$ *such that if $u$ is any $C^2(D)$ solution of*

$$-\Delta u + u|u|^q = 0 \qquad in \ D,$$

*then for all $x \in D$,*

$$|u(x)| \leq \frac{C}{\rho(x)^{2/q}}.$$

Recall that in our case $D$ is supposed to be bounded and $\partial D \in C^2$. The process $X^x$ is the solution of (4), and the stopping time $\tau_x$ is defined by (5). We will prove the following:

THEOREM 6 (A priori estimate). *There exists a constant $C$ [depending on the open set $D$, on $q$ and on the bound in (B) of $b$ and $\sigma$] such that for every $x \in \overline{D}$ and every solution $(Y, Z)$ of the BSDE (3) with terminal time $\tau_x$ and terminal data $\xi$ such that (H6) holds, we have*

$$\text{(12)} \qquad \forall t \geq 0 \qquad |Y_t| \leq \frac{C}{(\rho(X^x_{t \wedge \tau_x}))^{2/q}}.$$

We define the signed distance $d$

$$d(x) = \begin{cases} \text{dist}(x, \partial D) = \rho(x), & \text{if } x \in D, \\ -\text{dist}(x, \partial D), & \text{if } x \in \mathbb{R}^d \setminus D. \end{cases}$$

For $\mu > 0$, let

$$\Gamma_\mu \triangleq \{x \in \mathbb{R}^d; |d(x)| < \mu\}.$$

The following lemma (see [9], Lemma 14.16) relates the smoothness of the distance function $d$ in $\Gamma_\mu$ to that of the boundary $\partial D$.

LEMMA 1. *Let $D$ be bounded and $\partial D \in C^k$ for $k \geq 2$. Then there exists a positive constant $\mu$ depending on $D$ such that $d \in C^k(\Gamma_\mu)$.*

PROOF OF THEOREM 6. Recall that $D$ is an open bounded subset of $\mathbb{R}^d$ with $\partial D \in C^2$. From the previous lemma we already know that there exists a positive constant $\mu$ such that on $\Gamma_\mu$, the signed distance function $d$ belongs to $C^2$. And $d = \rho$ is continuous on $\overline{D}$. There exists a positive constant $R$ (depending only on $D$) such that for all $x \in \overline{D}$, $0 \leq d(x) = \rho(x) \leq R$. Let $\Phi \in C^\infty(\mathbb{R}^d; [0,1])$ such that $\Phi$ is equal to 1 on $\mathbb{R}^d \setminus \Gamma_\mu$ and is equal to 0 on $\Gamma_{\mu/2}$.

For $0 < \varepsilon \leq 1$ and $C > 0$, we define a function $\Psi_\varepsilon \in C^2(\mathbb{R}^d; \mathbb{R}_+)$ such that on $\overline{D}$,

$$\Psi_\varepsilon = \frac{C}{[(1 - \Phi)\rho + R\Phi + \varepsilon]^{2/q}}.$$

Such a function exists because $(1 - \Phi)\rho + R\Phi + \varepsilon \geq \varepsilon$ on $\overline{D}$. Remark that if $x \in \overline{D}$,

$$\Psi_\varepsilon(x) \leq \frac{C}{\rho(x)^{2/q}}.$$

We denote by $\theta_\varepsilon$ the function $(1 - \Phi)\rho + R\Phi + \varepsilon$, that is, on $\overline{D}$, $\Psi_\varepsilon = C\theta_\varepsilon^{-2/q}$. We apply the Itô formula to $\Psi_\varepsilon(X^x_{t \wedge \tau_x})$, where $x \in \overline{D}$. For convenience, we fix $\varepsilon > 0$ and $x \in \overline{D}$ and we omit the index $\varepsilon$ and $x$. For all

$0 \leq t \leq T$,

$$\Psi(X_{t \wedge \tau}) = \Psi(X_{T \wedge \tau}) - \int_{t \wedge \tau}^{T \wedge \tau} \Psi^{1+q}(X_r) \, dr$$

$$- \int_{t \wedge \tau}^{T \wedge \tau} \nabla \Psi(X_r) \sigma(X_r) \, dB_r$$

(13)

$$- \int_{t \wedge \tau}^{T \wedge \tau} [\nabla \Psi(X_r) b(X_r)$$

$$+ \tfrac{1}{2} \operatorname{Trace}(\sigma \sigma^*(X_r) D^2 \Psi(X_r)) - \Psi(X_r)^{1+q}] \, dr.$$

Now

$$\Psi^{1+q} = C^{1+q} \theta^{-2/q-2},$$

$$\frac{\partial \Psi}{\partial x_i} = -\frac{2C}{q} \theta^{-2/q-1} \frac{\partial \theta}{\partial x_i},$$

$$\frac{\partial^2 \Psi}{\partial x_i \partial x_j} = \frac{2C}{q}\left(\frac{2}{q}+1\right)\theta^{-2/q-2} \frac{\partial \theta}{\partial x_i} \frac{\partial \theta}{\partial x_j} - \frac{2C}{q} \theta^{-2/q-1} \frac{\partial^2 \theta}{\partial x_i \partial x_j}.$$

Therefore,

$$(\nabla \Psi) b + \frac{1}{2} \operatorname{Trace}(\sigma \sigma^* D^2 \Psi) - \Psi^{1+q}$$

(14)
$$= -C \theta^{-2/q-2} \bigg[ C^q + \frac{2\theta}{q}(\nabla \theta) b$$

$$- \frac{1}{q}\left(\frac{2}{q}+1\right) \|\sigma \nabla \theta\|^2 + \frac{\theta}{q} \operatorname{Trace}(\sigma \sigma^* D^2 \theta) \bigg];$$

and $b$, $\sigma$, $\theta$, $\nabla \theta$ and $D^2 \theta$ are bounded on $\overline{D}$. So we can choose the constant $C$ such that for all $x \in \overline{D}$

(15)
$$C^q + \frac{2\theta(x)}{q}(\nabla \theta(x)) b(x) - \frac{1}{q}\left(\frac{2}{q}+1\right) \|\sigma(x) \nabla \theta(x)\|^2$$

$$+ \frac{\theta(x)}{q} \operatorname{Trace}(\sigma(x) \sigma^*(x) D^2 \theta(x)) \geq 0.$$

The constant $C$ depends only on $D$, on $q$ and on the bound in (B) of $b$ and $\sigma$. We have obtained for all $0 \leq t \leq T$,

$$\Psi(X_{t \wedge \tau}) = \Psi(X_{T \wedge \tau}) - \int_{t \wedge \tau}^{T \wedge \tau} \nabla \Psi(X_r) \sigma(X_r) \, dB_r$$

$$- \int_{t \wedge \tau}^{T \wedge \tau} \Psi(X_r)^{1+q} \, dr + \int_{t \wedge \tau}^{T \wedge \tau} U_r \, dr;$$

with $U$ a nonnegative adapted process, and on $\{t \geq \tau\}$,

$$\Psi(X_{t\wedge\tau}) = \frac{C}{\varepsilon^{2/q}}.$$

If $(Y, Z)$ is the solution of the BSDE (3) with a final condition $\xi$ in $L^\infty(\Omega, \mathcal{F}_\tau, \mathbf{P})$ (see Remark 4), we can find $0 < \varepsilon < 1$ such that

$$|\xi| \leq \frac{C}{\varepsilon^{2/q}} \quad \text{a.s.}$$

Moreover, the Tanaka formula (see [10]) leads to, for all $0 \leq t \leq T$,

$$|Y_{t\wedge\tau}| = |Y_{T\wedge\tau}| - \int_{t\wedge\tau}^{T\wedge\tau} \operatorname{sign}(Y_r) Y_r |Y_r|^q \, dr - \int_{t\wedge\tau}^{T\wedge\tau} \operatorname{sign}(Y_r) Z_r \, dB_r$$
$$+ 2(\Lambda_{t\wedge\tau}(0) - \Lambda_{T\wedge\tau}(0))$$
$$= |Y_{T\wedge\tau}| - \int_{t\wedge\tau}^{T\wedge\tau} |Y_r|^{1+q} \, dr - \int_{t\wedge\tau}^{T\wedge\tau} \operatorname{sign}(Y_r) Z_r \, dB_r$$
$$+ 2(\Lambda_t(0) - \Lambda_T(0)),$$

where $\Lambda$ is a local time of $Y$. Thus, $\Lambda_{T\wedge\tau}(0) \geq \Lambda_{t\wedge\tau}(0)$ a.s.

By a comparison theorem (Corollary 4.4.2 in [6]), we have a.s.

$$\forall t \geq 0 \quad |Y_t| \leq \Psi_\varepsilon(X_{t\wedge\tau}) \leq \frac{C}{\rho(X_{t\wedge\tau})^{2/q}}.$$

By a density argument, it is clear that, if $(Y, Z)$ is the solution of the BSDE (3) with a final condition $\xi$ satisfying (H6), then

$$\forall t \geq 0 \quad |Y_t| \leq \frac{C}{\rho(X_{t\wedge\tau})^{2/q}}. \qquad \square$$

REMARK 3. *The constant $C$ depends only on the bound in* (B), *on $D$ and $q$. Moreover, in the special case where $D$ is a ball and where the drift in the SDE is equal to 0, the constant $C$ depends only on $q$ and on the bound of $\sigma$ and not on the center nor the radius of the ball. For example, if $X$ is the Brownian motion,* (12) *is true for any $C$ such that*

$$C^q \geq \frac{4}{q}\left(\frac{2}{q} + 1\right) + \frac{4d}{q}.$$

PROOF. In the case of a ball, we can give a slightly different proof because we have an explicit expression for the distance function. We will assume that $D$ is the ball centered at $y$ and with radius $R$. In this case the function $\rho$ is equal to

$$\rho(x) = R - |x - y| \leq \frac{R^2 - |x - y|^2}{R} = \theta(x),$$

if $x$ is in the ball. The function $\theta$ is not of class $C^2$ on the whole space. We modify $\theta$ in order to have a $C^2$ function. For all $0 < \varepsilon \leq R^2$, on the ball $\overline{G}$, $\theta_\varepsilon$ will be equal to

$$\theta_\varepsilon(x) = \frac{R^2 + \varepsilon - |x - y|^2}{R}$$

and on the whole space $\mathbb{R}^d$, $\theta_\varepsilon$ is positive and of class $C^2$. Remark that on $D$, $\theta_\varepsilon$ is greater than $\theta$. Now we consider the function

$$\Psi_\varepsilon(x) = \frac{C}{\theta_\varepsilon(x)^{2/q}}$$

and like in the proof of Theorem 6, we prove that there exists a constant $C$ such that the inequality (15) with $b=0$ holds. But now for $x \in \overline{D}$,

$$\nabla \theta_\varepsilon(x) = -\frac{2}{R}(x-y) \Rightarrow \|\sigma(x)\nabla\theta_\varepsilon(x)\|^2 = \frac{4}{R^2}\|\sigma(x)(x-y)\|^2$$

$$\leq \frac{4K^2}{R^2}|x-y|^2 \leq 4K^2;$$

$$D^2\theta_\varepsilon(x) = -\frac{2}{R}\mathrm{Id} \Rightarrow \mathrm{Trace}(\sigma(x)\sigma^*(x)D^2\theta_\varepsilon(x))$$

$$= -\frac{2}{R}\mathrm{Trace}(\sigma(x)\sigma^*(x)) \geq -\frac{2K}{R}.$$

It suffices to choose $C$ such that

$$C^q \geq \frac{4}{q}\left(\frac{2}{q}+1\right)K^2 + \frac{4K}{q}$$

in order to have (15). If $X$ is the Brownian motion $B$, that is, $\sigma = \mathrm{Id}$, $C^q \geq \frac{4}{q}(\frac{2}{q}+1) + \frac{4d}{q}$. Here $C$ depends only on the dimension $d$ and on $q$, like in the Keller–Osserman theorem. □

**2. Approximation.** We first prove a technical result which gives a sufficient condition on $\xi$ to insure existence and uniqueness of the solution of the BSDE (3). In our case for all $p > 1$, $\nu_p = 0$ and $f(t,0,0) = 0$.

PROPOSITION 2. *Under the condition* (C2) *on the first exit time $\tau$ of the diffusion $X$, let $\xi$ be an $\mathcal{F}_\tau$-measurable random variable such that $\xi \in L^r$ with $r > 2(1+q)$. Hence, there exists $p > 1$ and $\lambda > 0$ such that the condition* (H6) *is satisfied:*

(H6) $$\mathbb{E}\left[e^{p\lambda\tau}|\xi|^p + \int_0^\tau e^{p\lambda t}|\mathbb{E}(\xi|\mathcal{F}_t)|^{p(1+q)}\,dt\right] < +\infty.$$

PROOF. With $\alpha > 1$, $\gamma > 1$ such that $1/\alpha + 1/\gamma = 1$, the Hölder inequality leads to
$$\mathbb{E}(e^{p\lambda\tau}|\xi|^p) \leq [\mathbb{E}(e^{\alpha p\lambda\tau})]^{1/\alpha}[\mathbb{E}(|\xi|^{\gamma p})]^{1/\gamma}.$$

If $\xi \in L^r$ with $r > 1$, there exists $\gamma > 1$ and $p > 1$ such that $\gamma p \leq r$. From (C2), we can choose $\lambda > 0$ such that $\lambda\alpha p \leq \beta$.

For the rest of the condition (H6), we have $f(y) = -y|y|^p$ and thus,

$$\mathbb{E}\left[\int_0^\tau e^{p\lambda t}|f(\mathbb{E}^{\mathcal{F}_t}\xi)|^p \, dt\right] \leq \mathbb{E}\left[\int_0^\tau e^{p\lambda t}\mathbb{E}^{\mathcal{F}_t}(|\xi|^{p(1+q)}) \, dt\right]$$
$$\leq \mathbb{E}\left[\frac{e^{p\lambda\tau} - 1}{\lambda p} \sup_{t \in [0,\tau]} \mathbb{E}^{\mathcal{F}_t}(|\xi|^{p(1+q)})\right]$$
$$\leq \frac{1}{\lambda p}[\mathbb{E}e^{\alpha p\lambda\tau}]^{1/\alpha}\left[\mathbb{E} \sup_{t \in [0,\tau]} \mathbb{E}^{\mathcal{F}_t}(|\xi|^{\gamma p(1+q)})\right]^{1/\gamma}.$$

From the Burkholder–Davis–Gundy inequality we obtain
$$\mathbb{E} \sup_{t \in [0,\tau]} \mathbb{E}^{\mathcal{F}_t}(|\xi|^{\gamma p(1+q)}) \leq C[\mathbb{E}|\xi|^{2\gamma p(1+q)}]^{1/2}.$$

If $r > 2(1+q)$, we can choose $\gamma > 1$, $p > 1$ and $\lambda > 0$ sufficiently small such that $2\gamma p(1+q) < r$ and $\alpha\lambda p \leq \beta$. $\square$

REMARK 4. From the previous proposition, if $\xi \in L^r$ for some $r > 2(1+q)$, there exists $p > 1$ and $\lambda > 0$ such that the condition (H6) is satisfied. From Theorem 1, there exists a unique $L^p$-solution $(Y, Z)$ of the BSDE (3) which satisfies the estimate (2).

For terminal data $\xi$, if $(Y, Z)$ is a $L^p$-solution for some $p > 1$, then $(Y, Z)$ is also a $L^{p'}$-solution for all $1 < p' \leq p$. Therefore, if $(Y', Z')$ is a $L^{p'}$-solution for some $p' \leq p$, then we can easily prove that $(Y, Z) = (Y', Z')$.

If $\xi \in L^\infty$, from the proof of Proposition 2, (H6) holds for every $p > 1$ and $\lambda = \beta/p$.

We are interested in the case where $\xi$ is a nonnegative random variable with this new assumption:
$$\mathbf{P}(\xi = +\infty) > 0.$$
We still assume that the conditions (L), (B) and (E) hold, that $\tau = \tau_x$ is the exit time of the diffusion $X^x$ from the set $\overline{D}$, and that $\tau$ satisfies (C1) and (C2). From Remark 2, we can suppose $x$ to be in $D$ and for convenience, we omit the variable $x$. We suppose that $\xi$ is $\mathcal{F}_\tau$-measurable.

Now for each $n \in \mathbb{N}^*$, let $\xi_n = \xi \wedge n$ be our final condition. With Remark 4 we obtain the following:

LEMMA 2. *There exists a unique solution $(Y^n, Z^n)$ (in the sense of Definition 1) of the BSDE (3) with terminal data $\xi \wedge n$.*

From a comparison theorem (see Corollary 4.4.2 in [6]), we have a.s.

$$\forall t \geq 0, \ n \leq m \qquad 0 \leq Y_t^n \leq Y_t^m.$$

Define the progressively measurable process $Y$ by

$$(16) \qquad Y_t = \lim_{n \to \infty} Y_t^n \qquad \forall t \geq 0.$$

PROPOSITION 3. *The sequence $(Z^n)_{n \in \mathbb{N}^*}$ converges also to a process $Z$ and $(Y, Z)$ satisfies the assumptions* (D1) *and* (D2) *of Definition 2.*

PROOF. From (16), we already know that $(Y^n)_{n \in \mathbb{N}^*}$ converges to $Y$.

From the Itô formula and the Burkholder–Davis–Gundy inequality, there exists a constant $K$ such that for all $\eta > 0$, $n \geq m$ and $0 \leq s$,

$$\mathbb{E}\left(\sup_{t \in [0,s]} |Y_{t \wedge \tau_\eta}^n - Y_{t \wedge \tau_\eta}^m|^2\right) + \mathbb{E} \int_0^{s \wedge \tau_\eta} \|Z_r^n - Z_r^m\|^2 \, dr$$

$$\leq K\mathbb{E}(|Y_{s \wedge \tau_\eta}^n - Y_{s \wedge \tau_\eta}^m|^2).$$

See (8) for the definition of $\tau_\eta$. But with the inequality (12),

$$Y_{s \wedge \tau_\eta}^n \leq \frac{C}{\rho(X_{s \wedge \tau_\eta})^{2/q}} \leq \frac{C}{\eta^{2/q}},$$

and with the Lebesgue theorem, we conclude that $(Y_{\cdot \wedge \tau_\eta}^n, Z_{\cdot \wedge \tau_\eta}^n)_n$ converges to $(Y_{\cdot \wedge \tau_\eta}, Z_{\cdot \wedge \tau_\eta})$ in $L^2(\Omega; C(\mathbb{R}_+; \mathbb{R}_+)) \times L^2(\Omega \times \mathbb{R}^+)$ and $(Y_{\cdot \wedge \tau_\eta}^n)$ converges uniformly to $Y_{\cdot \wedge \tau_\eta}$.

Hence, $(Y, Z)$ satisfies the following equation: for all $0 \leq t \leq s$, for all $\eta > 0$,

$$Y_{t \wedge \tau_\eta} = Y_{s \wedge \tau_\eta} - \int_{t \wedge \tau_\eta}^{s \wedge \tau_\eta} (Y_r)^{1+q} \, dr - \int_{t \wedge \tau_\eta}^{s \wedge \tau_\eta} Z_r \, dB_r.$$

From this equation, with the Itô formula and the estimate (12), we deduce that

$$\mathbb{E}\left(\sup_{t \in [0,s]} |Y_{t \wedge \tau_\eta}|^2\right) + \mathbb{E} \int_0^{s \wedge \tau_\eta} \|Z_r\|^2 \, dr \leq K\mathbb{E}|Y_{s \wedge \tau_\eta}|^2$$

$$\leq \frac{K}{\rho^{2/q}(X_{s \wedge \tau_\eta})} \leq \frac{K}{\eta^{2/q}}.$$

Therefore, $(Y, Z)$ satisfies the conditions (D1) and (D2) of Definition 2. $\square$

From Definition 1, we also have that on the set $\{t \geq \tau\}$, $Y_t = \xi$ and $Z_t = 0$. With the monotonicity of the sequence $Y^n$, we can conclude that

$$\liminf_{t \to +\infty} Y_{t \wedge \tau} \geq \xi \quad \text{a.s.}$$

It remains to show the converse inequality,

$$\limsup_{t \to +\infty} Y_{t \wedge \tau} \leq \xi,$$

to have the last condition (D3). Without more assumptions on $\xi$, we cannot prove (D3). But we are able to give some other estimates on $Y$ and $Z$.

PROPOSITION 4. *For all $\varepsilon > 1$, there exists $K$ such that*

$$\mathbb{E} \int_0^\tau \|Z_r\|^2 \rho(X_r)^{4/q+\varepsilon} \, dr \leq K.$$

PROOF. We use again the notations $\Gamma_\mu$, the function $\theta = \theta_0 = (1-\Phi)\rho + R\Phi$ as in the proof of Theorem 6, and $\tau_\eta$ for $\eta < \mu$. Recall that $\theta \in C^2(\overline{D}; \mathbb{R})$ and on $\overline{D}$, $\theta = (1 - \Phi)\rho + R\Phi \geq d \geq 0$. Of course, $x \mapsto |\theta(x)|^{4/q+\varepsilon}$ is not in $C^2(\mathbb{R}^d)$, but this function belongs to $C^2(D \setminus \Gamma_\eta)$ and we can define this function on the rest of $(\mathbb{R}^d \setminus D) \cup \Gamma_\eta$ in order to have the required regularity. The Itô formula leads to

$$(Y_{t \wedge \tau_\eta}^n)^2 \theta(X_{t \wedge \tau_\eta})^{4/q+\varepsilon}$$

$$= (Y_0^n)^2 \theta(X_0)^{4/q+\varepsilon} + \int_0^{t \wedge \tau_\eta} \|Z_r^n\|^2 \theta(X_r)^{4/q+\varepsilon} \, dr$$

$$+ 2 \int_0^{t \wedge \tau_\eta} (Y_r^n)^{2+q} \theta(X_r)^{4/q+\varepsilon} \, dr + 2 \int_0^{t \wedge \tau_\eta} Y_r^n \theta(X_r)^{4/q+\varepsilon} Z_r^n \, dB_r$$

$$+ \left(\frac{4}{q} + \varepsilon\right) \int_0^{t \wedge \tau_\eta} (Y_r^n)^2 \theta(X_r)^{4/q+\varepsilon-1} \nabla \theta(X_r) (b(X_r) \, dr + \sigma(X_r) \, dB_r)$$

$$+ \frac{(4/q+\varepsilon)}{2} \int_0^{t \wedge \tau_\eta} (Y_r^n)^2 \left[ \left(\frac{4}{q} + \varepsilon - 1\right) \theta(X_r)^{4/q+\varepsilon-2} \|\sigma(X_r) \nabla \theta(X_r)\|^2 \right.$$

$$\left. + \theta(X_r)^{4/q+\varepsilon-1} \operatorname{Trace}(\sigma \sigma^*(X_r) D^2 f(X_r)) \right] dr$$

$$+ 2 \left(\frac{4}{q} + \varepsilon\right) \int_0^{t \wedge \tau_\eta} Y_r^n \theta(X_r)^{4/q+\varepsilon-1} Z_r^n \nabla \theta(X_r) \sigma(X_r) \, dr.$$

Then

$$\mathbb{E} \int_0^{t \wedge \tau_\eta} \|Z_r^n\|^2 \theta(X_r)^{4/q+\varepsilon} \, dr$$

$$+ 2\left(\frac{4}{q} + \varepsilon\right) \mathbb{E} \int_0^{t \wedge \tau_\eta} Y_r^n \theta(X_r)^{4/q+\varepsilon-1} Z_r^n \nabla \theta(X_r) \sigma(X_r) \, dr$$

is bounded from above by

$$\mathbb{E}((Y_{t\wedge\tau_\eta}^n)^2\theta(X_{t\wedge\tau_\eta})^{4/q+\varepsilon})$$

$$-\left(\frac{4}{q}+\varepsilon\right)\mathbb{E}\int_0^{t\wedge\tau_\eta}(Y_r^n)^2\theta(X_r)^{4/q+\varepsilon-1}\nabla\theta(X_r)b(X_r)\,dr$$

(17) $$-\frac{1}{2}\left(\frac{4}{q}+\varepsilon\right)\mathbb{E}\int_0^{t\wedge\tau_\eta}(Y_r^n)^2\theta(X_r)^{4/q+\varepsilon-1}\operatorname{Trace}(\sigma\sigma^*(X_r)D^2\theta(X_r))\,dr$$

$$-\frac{1}{2}\left(\frac{4}{q}+\varepsilon\right)\left(\frac{4}{q}+\varepsilon-1\right)\mathbb{E}\int_0^{t\wedge\tau_\eta}(Y_r^n)^2 f(X_r)^{4/q}\theta(X_r)^{\varepsilon-2}$$

$$\times\|\sigma(X_r)\nabla\theta(X_r)\|^2\,dr.$$

In the proof of Theorem 6, we have obtained that there exists some constant $C$ such that for all $n\in\mathbb{N}^*$ and for all $t\geq 0$,

(18) $$(Y_{t\wedge\tau}^n)^2\theta(X_{t\wedge\tau})^{4/q}\leq C.$$

Moreover $b$ and $\sigma$ are bounded [assumption (B)], and $\nabla\theta$ and $D^2\theta$ are also bounded on $\overline{D}$. Thus, the right-hand side of (17) is bounded by

$$K\left(1+\mathbb{E}\int_0^\tau\theta^{\varepsilon-1}(X_r)\,dr+\mathbb{E}\int_0^\tau\theta^{\varepsilon-2}(X_r)\,dr\right).$$

We denote by $p(t,x,y)$ the density of $\mathbf{P}^x(X_t\in dy;\tau>t)$. $\mathbf{P}^x$ means that the diffusion process $X$ starts from $x\in D$ at time 0. Then

$$\mathbb{E}\int_0^\tau\theta^{\varepsilon-1}(X_r)\,dr=\int_0^\infty\int_D\theta^{\varepsilon-1}(y)p(r,x,y)\,dy\,dr$$

$$=\int_D\theta^{\varepsilon-1}(y)G(x,y)\,dy,$$

where $G$ is the Green function associated to the process $X$ killed at time $\tau$ (see [21], Section 4.2, Theorem 2.5).

We claim that the last integral is finite. Indeed, if $B(x,\nu)$ is the ball centered at $x$ with radius $\nu>0$, since $x\in D$, we can find $\nu>0$ such that $B(x,\nu)\subset D$ and $B(x,\nu)\cap\Gamma_\nu=\varnothing$. We denote by $U$ the set $D\setminus(B(x,\nu)\cup\Gamma_\nu)$. On $U$, $f^{\varepsilon-1}G(x,\cdot)$ is a continuous function and is bounded by $K$. On $\Gamma_\nu$ [resp. on $B(x,\nu)$], $G(x,\cdot)$ (resp. $\theta^{\varepsilon-1}$) is continuous and bounded by $K$. On the boundary of $D$, $\theta^{\varepsilon-1}$ is singular if $\varepsilon<1$. Recall the definition of $\theta:\theta=(1-\varphi)\rho+R\varphi$. Hence, $\theta$ is equivalent to $\rho$ at the boundary and if $\varepsilon>0$, $f^{\varepsilon-1}$ is integrable on $D$. $G$ has a singularity when $y\to x$, but with Theorem 2.8 and Exercise 4.16 in [21], this singularity is integrable. Therefore, we split the integral into three terms:

$$\mathbb{E}\int_0^\tau\theta^{\varepsilon-1}(X_r)\,dr=\int_D\theta^{\varepsilon-1}(y)G(x,y)\,dy$$

$$\leq \int_{B(x,\nu)} \theta^{\varepsilon-1}(y) G(x,y)\, dy + \int_{\Gamma_\nu} \theta^{\varepsilon-1}(y) G(x,y)\, dy$$

$$+ \int_U \theta^{\varepsilon-1}(y) G(x,y)\, dy$$

$$\leq K \int_{B(x,\nu)} G(x,y)\, dy + K \int_{\Gamma_\nu} \theta^{\varepsilon-1}(y)\, dy + K \operatorname{Vol}(D)$$

$$< +\infty.$$

From the second integral, the same arguments show that if $\varepsilon > 1$,

$$\mathbb{E} \int_0^\tau \theta^{\varepsilon-2}(X_r)\, dr < +\infty. \tag{19}$$

Therefore, the right-hand side of (17) is bounded by a constant $K$ which does not depend on $\eta$, $n$ and $t$. And using the Cauchy–Schwarz inequality,

$$\left| \mathbb{E} \int_0^{t \wedge \tau_\eta} Y_r^n \theta(X_r)^{4/q+\varepsilon-1} Z_r^n \nabla \theta(X_r) \sigma(X_r)\, dr \right|$$

$$\leq \left( \mathbb{E} \int_0^{t \wedge \tau_\eta} \|Z_r^n\|^2 \theta(X_r)^{4/q+\varepsilon}\, dr \right)^{1/2}$$

$$\times \left( \mathbb{E} \int_0^{t \wedge \tau_\eta} (Y_r^n)^2 \theta(X_r)^{4/q+\varepsilon-2} \|\nabla \theta(X_r) \sigma(X_r)\|^2\, dr \right)^{1/2}.$$

But since $\nabla \theta$ and $\sigma$ are bounded and since (18) and (19) hold, if $\varepsilon > 1$,

$$\mathbb{E} \int_0^{t \wedge \tau_\eta} (Y_r^n)^2 \theta(X_r)^{4/q+\varepsilon-2} \|\nabla \theta(X_r) \sigma(X_r)\|^2\, dr \leq K.$$

Inequality (17) can be written as follows: $A_n + B_n \leq C_n$ with

$$0 \leq A_n = \mathbb{E} \int_0^{t \wedge \tau_\eta} \|Z_r^n\|^2 \theta(X_r)^{4/q+\varepsilon}\, dr,$$

$$|B_n| = 2\left(\frac{4}{q} + \varepsilon\right) \left| \mathbb{E} \int_0^{t \wedge \tau_\eta} Y_r^n \theta(X_r)^{4/q+\varepsilon-1} Z_r^n \nabla \theta(X_r) \sigma(X_r)\, dr \right| \leq K A_n^{1/2}$$

and $|C_n| \leq K$. Thus, for all $n \in \mathbb{N}^*$ and for all $t \geq 0$,

$$\mathbb{E} \int_0^{t \wedge \tau_\eta} \|Z_r^n\|^2 \theta(X_r)^{4/q+\varepsilon}\, dr \leq K,$$

which implies, by Fatou's lemma,

$$\mathbb{E} \int_0^\tau \|Z_r\|^2 \theta(X_r)^{4/q+\varepsilon}\, dr \leq K.$$

Since $\theta \geq \rho$ on $\overline{D}$, we obtain the announced result. $\square$

The condition $\varepsilon > 1$ is required in order to insure that
$$\mathbb{E}\int_0^\tau \theta^{\varepsilon-2}(X_r)\,dr < \infty.$$
But this integral is equal to $\int_D \theta^{\varepsilon-2}(y)G(x,y)\,dy$, where $G(x,y)$ is the Green function associated with the process $X^x$ killed at $\tau$. And if, for example, the infinitesimal generator of the diffusion $X$ is self-adjoint in $L^2(\mathbb{R}^d)$, that is, $\mathcal{L} = (1/2)\operatorname{div}(\sigma\sigma^*\nabla)$, then $G(x,y) \leq K\rho(y)$ (see [7], Theorem 9.5) and the previous integral is finite for any $\varepsilon > 0$. Hence, we obtain that, for any $\varepsilon > 0$,
$$\mathbb{E}\int_0^\tau \|Z_r\|^2 \rho(X_r)^{4/q+\varepsilon}\,dr < \infty.$$

In the next proposition we find an adapted process smaller than $Y$. This process will give us a lower bound on the explosion rate of $Y$ on the blow-up set $\{\xi = \infty\}$.

PROPOSITION 5 (Lower bound on $Y$). *We define the following process:*
$$\Xi_t = \mathbb{E}^{\mathcal{F}_t}\left[\left(\frac{1}{q(\tau - \tau\wedge t) + 1/\xi^q}\right)^{1/q}\right].$$
*Then for all $t \geq 0$, $\Xi_t \leq Y_t$.*

PROOF. Denote by $\alpha_t$ the quantity
$$\alpha_t = \left(\frac{1}{q(\tau - \tau\wedge t) + 1/\xi^q}\right)^{1/q},$$
if $t < \tau$ and $\alpha_t = \xi$ on $\{t \geq \tau\}$. The process $\alpha$ solves the equation
$$\alpha_t = \alpha_T - \int_t^T \alpha_r^{1+q}\mathbf{1}_{[0,\tau]}(r)\,dr.$$
Note that $\Xi_t = \mathbb{E}(\alpha_t|\mathcal{F}_t)$. Thanks to Jensen's inequality,
$$\Xi_t \leq \mathbb{E}^{\mathcal{F}_t}\left(\Xi_T - \int_t^T \Xi_r^{1+q}\mathbf{1}_{[0,\tau]}(r)\,dr\right)$$
and the comparison theorem (Corollary 4.4.2 in [6]) achieves the proof. □

We now want to find a lower bound for $\rho(X_{t\wedge\tau})^{2/q}\Xi_t$ when $t$ goes to $\infty$ on $\{\xi = \infty\}$.

LEMMA 3. *Let $\rho(x)$ denote the distance of $x \in D$ to the boundary $\partial D$ and $\tau$ be the exit time from $\overline{D}$ of the diffusion $X$. Let the conditions* (L), (B) *and* (E) *hold. Then there exist two positive constants $C_1$ and $C_2$ which depend on $D$, $q$, $\sigma$ and $b$ such that for all $x \in D$,*
$$C_1 \leq \rho(x)^{2/q}\mathbb{E}_x\left[\left(\frac{1}{\tau}\right)^{1/q}\right] \leq C_2.$$

PROOF. Recall that if $x \in D$, $\mathbf{P}^x(\tau > 0) = 1$ and
$$\mathbb{E}_x\left(\frac{1}{\tau^{1/q}}\right) = \int_0^{+\infty} \mathbf{P}^x\left(\tau < \frac{1}{y^q}\right) dy.$$
If $\tau < h$, then $\sup_{t \in [0,h]} |X_t - x| > \rho(x)$. Therefore, we can apply Theorem 4.2.1, page 87 of [23] to obtain
$$\mathbf{P}^x(\tau < h) \leq \mathbf{P}^x\left(\sup_{t \in [0,h]} |X_t - x| > \rho(x)\right) \leq K_1 e^{K_2 h} e^{-K_2 \rho(x)^2/h}.$$
We apply this inequality with $h = 1/y^q$ and $y \geq 1$:
$$\mathbb{E}_x\left(\frac{1}{\tau^{1/q}}\right) \leq 1 + \int_1^{+\infty} K_1 e^{K_2/y^q} e^{-K_2 \rho(x)^2 y^q} dy$$
$$\leq 1 + K_1 e^{K_2} \int_1^{+\infty} e^{-K_2 \rho(x)^2 y^q} dy$$
$$= 1 + \frac{K_1 e^{K_2}}{q\rho(x)^{2/q}} \int_{\rho(x)^2}^{+\infty} e^{-K_2 u} u^{1/q-1} du$$
$$\leq 1 + \frac{K_1 e^{K_2}}{q\rho(x)^{2/q}} \int_0^{+\infty} e^{-K_2 u} u^{1/q-1} du.$$
Since $-1 + 1/q > -1$, we deduce
$$\rho(x)^{2/q} \mathbb{E}_x\left(\frac{1}{\tau^{1/q}}\right) \leq C_2.$$

For the other inequality remark that
$$\int_0^{+\infty} \mathbf{P}^x\left(\tau < \frac{1}{y^q}\right) dy \geq \int_0^{1/\rho(x)^{2/q}} \mathbf{P}^x\left(\tau < \frac{1}{y^q}\right) dy$$
and $\mathbf{P}^x(\tau < 1/y^q) \geq \mathbf{P}^x(X_{1/y^q} \notin \overline{D})$. We just have to find a lower bound to
$$\rho(x)^{2/q} \int_0^{1/\rho(x)^{2/q}} \mathbf{P}^x(X_{1/y^q} \notin \overline{D}) dy = \int_1^{+\infty} \mathbf{P}^x(X_{\rho(x)^2 u^q} \notin \overline{D}) \frac{du}{u^2}.$$
Let $\Delta_x$ be the set $(\mathbb{R}^d \setminus \overline{D}) \cap B(x, 2\rho(x))$ which is not empty, and $\text{Vol}(\Delta_x)$ denotes the volume of the set $\Delta_x$. Using the Aronson estimates of [22], we have, for $u \geq 1$,
$$\mathbf{P}^x(X_{\rho(x)^2 u^q} \notin \overline{D})$$
$$\geq K_3 \int_{\mathbb{R}^d \setminus \overline{D}} \left(\frac{1}{2\pi u^q \rho(x)^2}\right)^{d/2}$$
$$\times \exp\left(-K_4 u^q \rho(x)^2 - K_4 \frac{|y-x|^2}{2u^q \rho(x)^2}\right) dy$$

$$\geq K_3 \int_{\Delta_x} \left(\frac{1}{2\pi u^q \rho(x)^2}\right)^{d/2}$$
$$\times \exp\left(-K_4 u^q \rho(x)^2 - K_4 \frac{|y-x|^2}{2u^q \rho(x)^2}\right) dy$$
$$\geq K_3 \left(\frac{1}{2\pi u^q \rho(x)^2}\right)^{d/2} \exp(-K_4 u^q \rho(x)^2) \int_{\Delta_x} \exp(-2K_4 u^{-q}) \, dy$$
$$\geq K_3 e^{-2K_4} \left(\frac{1}{2\pi u^q}\right)^{d/2} \exp(-K_4 u^q \rho(x)^2) \frac{\mathrm{Vol}(\Delta_x)}{\rho(x)^d}$$

because $u \geq 1$. Thus,

$$\rho(x)^{2/q} \int_0^{1/\rho(x)^{2/q}} \mathbf{P}^x(X_{1/y^q} \notin \overline{D}) \, dy$$
$$\geq K_5 \left[\int_1^{+\infty} \exp(-K_4 u^q \rho(x)^2) \frac{du}{u^{2+dq/2}}\right] \frac{\mathrm{Vol}(\Delta_x)}{\rho(x)^d}$$

with
$$K_5 = K_3 e^{-2K_4} \left(\frac{1}{2\pi}\right)^{d/2}.$$

The integral has a lower bound because the open set $D$ is bounded and the dominated convergence theorem shows that

$$\lim_{\rho(x) \to 0} \int_1^{+\infty} \exp(-K_4 u \rho(x)^2) \frac{du}{u^{2+d/2}} = \int_1^{+\infty} \frac{du}{u^{2+d/2}} = K_6.$$

We have supposed that $\partial D \in C^2$ (which was important in the proof of Theorem 6). Therefore, the curvature is continuous on $\partial D$ which is compact; so the curvature is bounded. There exists $r > 0$ such that each point $y \in \partial D$ lies on the boundary of a ball with radius $r$ and this ball is contained in the complementary of $D$ (see Figure 1).

Instead of calculating the volume of $(\mathbb{R}^d \setminus \overline{D}) \cap B(x, 2\rho(x))$, our problem is reduced to the following: we find the volume of the intersection of two balls in that case; see Figure 2.

If $r < x < 3r$, the volume is equal to

$$V(x) = Cr^d \int_0^{\alpha(x)} (\sin\theta)^d \, d\theta + C2^d (x-r)^d \int_0^{\beta(x)} (\sin\theta)^d \, d\theta,$$

where $C$ is the volume of the unit ball in $\mathbb{R}^{d-1}$,

$$\alpha(x) = \mathrm{Arccos}\left(\frac{x^2 + r^2 - 4(x-r)^2}{2xr}\right)$$

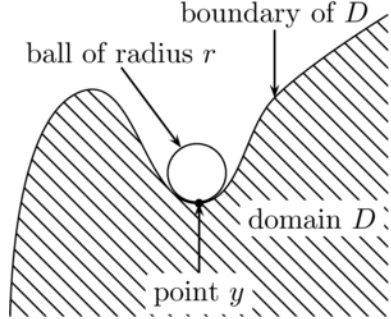

Fig. 1.

and
$$\beta(x) = \text{Arccos}\left(\frac{5x - 3r}{4x}\right).$$

Now we must prove that $\frac{V(x)}{(x-r)^d} \geq K_1$. We split $V(x)$ in two parts. For the first part, $C2^d(x-r)^d \int_0^{\beta(x)} (\sin\theta)^d \, d\theta$, the result is clear because if $r \leq x \leq 2r$,
$$0 < \text{Arccos}(7/8) \leq \beta(x) \leq \frac{\pi}{3}.$$

For the second part, $Cr^d \int_0^{\alpha(x)} (\sin\theta)^d \, d\theta$, if $r \leq x \leq \frac{4+\sqrt{7}}{3}r$,
$$\frac{x^2 + r^2 - 4(x-r)^2}{2xr} \in [0, 1] \implies \alpha(x) \in [0, \pi/2]$$
and we use the fact that sin is an increasing function on $[0, \pi/2[$; so
$$\int_0^{\alpha(x)} (\sin\theta)^d \, d\theta \leq \alpha(x)\left(1 - \left(\frac{x^2 + r^2 - 4(x-r)^2}{2xr}\right)^2\right)^{d/2}$$

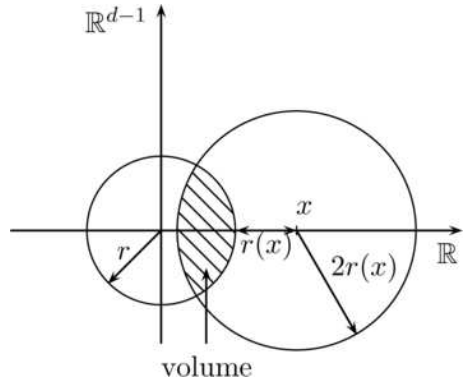

Fig. 2.

$$\leq 2^{d/2}\alpha(x)\left(1 - \frac{x^2 + r^2 - 4(x-r)^2}{2xr}\right)^{d/2}$$

$$= \left(\frac{3}{xr}\right)^{d/2} \alpha(x)(x-r)^d.$$

Therefore, if $r < x \leq \frac{4+\sqrt{7}}{3}r$,

$$\frac{V(x)}{(x-r)^d} \geq C2^d \int_0^{\text{Arccos}(7/8)} (\sin\theta)^d \, d\theta > 0.$$

This proves that $\frac{V(x)}{(x-r)^d} \geq \widetilde{K}$ and therefore, $\rho(x)^2 \mathbb{E}_x(\frac{1}{\tau}) \geq C_1 = K_5 \widetilde{K}$. This achieves the proof in the uniformly elliptic case. $\square$

REMARK 5. If the diffusion matrix is degenerate, the result on the lower bound may be false. Suppose that $\sigma \equiv 0$ and $b$ is bounded by $k$. If the exit time $\tau$ is smaller than $1/y^q$,

$$\frac{k}{y^q} \geq \int_0^{1/y^q} |b(X_r)| \, dr \geq \sup_{[0,1/y^q]} \left|\int_0^t b(X_r) \, dr\right| = \sup_{[0,1/y^q]} |X_t - x| > \rho(x)$$

and thus,

$$\rho(x)^{2/q} \mathbb{E}_x\left(\frac{1}{\tau^{1/q}}\right) \leq k\rho(x)^{1/q}$$

and the limit, as $\rho(x)$ goes to zero, is zero.

From the inequality (12), we already know that there exists a constant $C$ such that

$$\forall t \geq 0 \qquad Y_t \leq \frac{C}{(\rho(X_{t \wedge \tau}))^{2/q}}.$$

Now we prove Proposition 1:

PROOF OF PROPOSITION 1. From Proposition 5, we work on the process $\Xi = (\Xi_t)_{t \geq 0}$. For all $t \geq 0$,

$$\rho^{2/q}(X_{t\wedge\tau})\Xi_t = \rho^{2/q}(X_{t\wedge\tau})\mathbb{E}^{\mathcal{F}_t}\left[\left(\frac{1}{q(\tau - \tau \wedge t) + 1/\xi^q}\right)^{1/q}\right]$$

$$= \rho^{2/q}(X_{t\wedge\tau})\mathbb{E}^{\mathcal{F}_t}\left[\left(\frac{\xi^q \mathbf{1}_{\xi<\infty}}{1 + q\xi^q(\tau - \tau \wedge t)}\right)^{1/q}\right]$$

$$+ \rho^{2/q}(X_{t\wedge\tau})\mathbb{E}^{\mathcal{F}_t}\left[\left(\frac{1}{q(\tau - \tau \wedge t)}\right)^{1/q} \mathbf{1}_{\xi=\infty}\right].$$

The first term in the right-hand side is nonnegative. Let $\widetilde{\tau} = \tau - \tau \wedge t$: it is the first exit time of the diffusion starting at $X_{t \wedge \tau}$. Hence,

$$\rho^{2/q}(X_{t\wedge\tau})\mathbb{E}^{\mathcal{F}_t}\left[\left(\frac{1}{q(\tau - \tau \wedge t)}\right)^{1/q}\mathbf{1}_{\xi=\infty}\right]$$
$$= \left(\frac{1}{q}\right)^{1/q}\mathbb{E}^{\mathcal{F}_t}\left\{\rho^{2/q}(X_{t\wedge\tau})\mathbb{E}^{X_{t\wedge\tau}}\left[\left(\frac{1}{\widetilde{\tau}}\right)^{1/q}\right]\mathbf{1}_{\xi=\infty}\right\} \geq C_1\left(\frac{1}{q}\right)^{1/q}\mathbb{E}^{\mathcal{F}_t}(\mathbf{1}_{\xi=\infty}),$$

where $C_1$ is the lower bound of Lemma 3. Thus, we obtain

$$\rho^{2/q}(X_{t\wedge\tau})\Xi_t \geq C_1\left(\frac{1}{q}\right)^{1/q}\mathbb{E}^{\mathcal{F}_t}(\mathbf{1}_{\xi=\infty})$$

and we deduce the announced result. □

**3. Continuity.** Recall that we have constructed a couple of processes $(Y, Z)$ which satisfy for all $\eta > 0$ and all $0 \leq t \leq T$, $Y_t \geq 0$ and

$$Y_{t \wedge \tau_\eta} = Y_{T \wedge \tau_\eta} - \int_{t\wedge\tau_\eta}^{T\wedge\tau_\eta} (Y_r)^{1+q} \, dr - \int_{t\wedge\tau_\eta}^{T\wedge\tau_\eta} Z_r \, dB_r.$$

Moreover, on the set $\{t \geq \tau\}$, $Y_t = \xi$, $Z_t = 0$ and $\liminf_{t \to +\infty} Y_{t \wedge \tau} \geq \xi$ a.s. We now want to prove the converse inequality, namely, $\limsup_{t \to +\infty} Y_{t \wedge \tau} \leq \xi$ a.s. Remark that we just have to show this estimate on the set $\{\xi < +\infty\}$.

3.1. *Existence of the limit.* We first prove that the limit of $Y_{t \wedge \tau}$, as $t$ goes to $+\infty$, exists a.s. In the proof we will distinguish the two cases: $\xi$ is greater than a positive constant and $\xi$ is nonnegative.

3.1.1. *The case where $\xi$ is bounded away from zero.* We can show that $(Y_{t\wedge\tau})_{t\geq 0}$ has a limit as $t \to +\infty$ by using Itô's formula applied to the process $1/(Y^n)^q$. We prove the following result:

PROPOSITION 6. *Let the conditions* (B) *and* (E) *hold. Suppose there exists a real $\alpha > 0$ such that $\xi \geq \alpha > 0$, $\mathbf{P}$-a.s. Then*

$$(20) \qquad Y_{t\wedge\tau} = \left[\mathbb{E}^{\mathcal{F}_t}\left(q(\tau - t \wedge \tau) + \left(\frac{1}{\xi^q}\right)\right) - \Phi_t\right]^{-1/q}, \qquad 0 \leq t,$$

*where $\Phi$ is a nonnegative supermartingale such that on the set $\{t \geq \tau\}$, $\Phi_t = 0$.*

PROOF. From Proposition 5, for every $n \in \mathbb{N}^*$ and every $0 \leq t$,

$$Y_t^n \geq \Xi_t^n = \mathbb{E}^{\mathcal{F}_t}\left[\left(\frac{1}{q(\tau - \tau \wedge t) + (1/(\xi \wedge n))^q}\right)^{1/q}\right].$$

Since $\xi \geq \alpha$, we have

$$\Xi_t^n \geq \mathbb{E}^{\mathcal{F}_t}\left[\left(\frac{1}{q(\tau - \tau \wedge t) + (1/\alpha)^q}\right)^{1/q}\right] \geq \alpha \mathbb{E}^{\mathcal{F}_t}\left[\left(\frac{1}{1 + q\tau\alpha^q}\right)^{1/q}\right]$$

$$\geq \alpha\left(\frac{1}{1 + q\alpha^q \mathbb{E}^{\mathcal{F}_t}(\tau)}\right)^{1/q}.$$

Therefore,

$$(21) \qquad \forall t \geq 0 \qquad 0 \leq \frac{1}{(Y_t^n)^q} \leq \frac{1}{\alpha^q}(1 + q\alpha^q \mathbb{E}^{\mathcal{F}_t}(\tau)) < +\infty$$

because the conditions (B) and (E) hold, which implies, in particular, that $\tau \in L^1(\Omega)$. Thus, for all $t \geq 0$, $(Y_t^n)^{-q}$ belongs to $L^1(\Omega)$. We want to apply the Itô formula to the semi-martingale $Y^n$ with the function $0 < x \mapsto x^{-q}$. But we just have that a.s. for all $t \geq 0$, $Y_t^n > 0$. For $\varepsilon > 0$, we define a $C^2$ function $f_\varepsilon : \mathbb{R} \to \mathbb{R}$ such that on $\mathbb{R}_+$,

$$f_\varepsilon(x) = \left(\frac{1}{x + \varepsilon}\right)^q.$$

Note that for a fixed $x \in \mathbb{R}_+$, $(f_\varepsilon(x))_{\varepsilon > 0}$ is increasing and the limit is equal to $f(x) = x^{-q}$. By the Itô formula, for all $0 \leq t \leq T$,

$$f_\varepsilon(Y_{t \wedge \tau}^n) = f_\varepsilon(Y_{T \wedge \tau}^n) - \int_{t \wedge \tau}^{T \wedge \tau} f_\varepsilon'(Y_r^n) Z_r^n \, dB_r - \int_{t \wedge \tau}^{T \wedge \tau} f_\varepsilon'(Y_r^n)(Y_r^n)^{1+q} \, dr$$

$$- \tfrac{1}{2} \int_{t \wedge \tau}^{T \wedge \tau} f_\varepsilon''(Y_r^n) \|Z_r^n\|^2 \, dr$$

$$(22)$$

$$= \mathbb{E}^{\mathcal{F}_t} f_\varepsilon(Y_{T \wedge \tau}^n) - \mathbb{E}^{\mathcal{F}_t} \int_{t \wedge \tau}^{T \wedge \tau} f_\varepsilon'(Y_r^n)(Y_r^n)^{1+q} \, dr$$

$$- \tfrac{1}{2} \mathbb{E}^{\mathcal{F}_t} \int_{t \wedge \tau}^{T \wedge \tau} f_\varepsilon''(Y_r^n) \|Z_r^n\|^2 \, dr.$$

Now for $x \geq 0$, $f_\varepsilon'(x) x^{1+q} = -q(\frac{x}{x+\varepsilon})^{1+q}$, thus, $-q \leq f_\varepsilon'(x) x^{1+q} \leq 0$, and

$$f_\varepsilon''(x) = q(1 + q)\left(\frac{1}{x + \varepsilon}\right)^{2+q}.$$

Thereby, a.s. and in $L^1(\Omega)$ for all $0 \leq t \leq T$,

$$\lim_{\varepsilon \to 0} \mathbb{E}^{\mathcal{F}_t} \int_{t \wedge \tau}^{T \wedge \tau} f_\varepsilon'(Y_r^n)(Y_r^n)^{1+q} \, dr = -q \mathbb{E}^{\mathcal{F}_t}(T \wedge \tau - t \wedge \tau).$$

From (21), we have that a.s. and in $L^1(\Omega)$

$$\lim_{\varepsilon \to 0} \mathbb{E}^{\mathcal{F}_t} f_\varepsilon(Y_{T \wedge \tau}^n) = \mathbb{E}^{\mathcal{F}_t} \frac{1}{(Y_{T \wedge \tau}^n)^q}$$

and

$$\lim_{\varepsilon \to 0} f_\varepsilon(Y^n_{t \wedge \tau}) = \frac{1}{(Y^n_{t \wedge \tau})^q}.$$

For the last term in (22), we use the monotone convergence theorem and hence, we have proved that, for all $0 \leq t \leq T$,

$$\frac{1}{(Y^n_{t\wedge\tau})^q} = \mathbb{E}^{\mathcal{F}_t} \frac{1}{(Y^n_{T\wedge\tau})^q} + q\mathbb{E}^{\mathcal{F}_t}(T \wedge \tau - t \wedge \tau)$$
$$- \frac{q(q+1)}{2} \mathbb{E}^{\mathcal{F}_t} \int_{t\wedge\tau}^{T\wedge\tau} \frac{\|Z^n_r\|^2}{(Y^n_r)^{2+q}} dr.$$

Let $T$ go to $+\infty$:

(23)
$$\frac{1}{(Y^n_{t\wedge\tau})^q} = \mathbb{E}^{\mathcal{F}_t} \frac{1}{(\xi \wedge n)^q} + q\mathbb{E}^{\mathcal{F}_t}(\tau - t \wedge \tau)$$
$$- \frac{q(q+1)}{2} \mathbb{E}^{\mathcal{F}_t} \int_{t\wedge\tau}^{\tau} \frac{\|Z^n_r\|^2}{(Y^n_r)^{2+q}} dr.$$

Let $n \geq m$. Since $\xi \wedge n \geq \xi \wedge m$, we obtain, for all $0 \leq t$,

$$0 \leq \frac{1}{(Y^m_{t\wedge\tau})^q} - \frac{1}{(Y^n_{t\wedge\tau})^q}$$
$$= \mathbb{E}^{\mathcal{F}_t}\left(\frac{1}{(\xi \wedge m)^q} - \frac{1}{(\xi \wedge n)^q}\right)$$
$$- \frac{q(q+1)}{2}\left(\mathbb{E}^{\mathcal{F}_t} \int_{t\wedge\tau}^{\tau} \frac{\|Z^m_s\|^2}{(Y^m_s)^{q+2}} ds - \mathbb{E}^{\mathcal{F}_t} \int_{t\wedge\tau}^{\tau} \frac{\|Z^n_s\|^2}{(Y^n_s)^{q+2}} ds\right).$$

Now

$$\frac{q(q+1)}{2}\left|\mathbb{E}^{\mathcal{F}_t} \int_{t\wedge\tau}^{\tau} \frac{\|Z^m_s\|^2}{(Y^m_s)^{q+2}} ds - \mathbb{E}^{\mathcal{F}_t} \int_{t\wedge\tau}^{\tau} \frac{\|Z^n_s\|^2}{(Y^n_s)^{q+2}} ds\right|$$
$$\leq \left[\mathbb{E}^{\mathcal{F}_t}\left(\frac{1}{(\xi \wedge m)^q} - \frac{1}{(\xi \wedge n)^q}\right)\right] \vee \left[\frac{1}{(Y^m_{t\wedge\tau})^q} - \frac{1}{(Y^n_{t\wedge\tau})^q}\right].$$

For a fixed $t \geq 0$, the sequences $(\mathbb{E}^{\mathcal{F}_t} \frac{1}{(\xi \wedge n)^q})_{n \geq 1}$ and $(\frac{1}{(Y^n_{t\wedge\tau})^q})_{n \geq 1}$ converge a.s. and in $L^1$ (dominated convergence theorem). Then the sequence $(\mathbb{E}^{\mathcal{F}_t} \int_{t\wedge\tau}^{\tau} \frac{\|Z^n_s\|^2}{(Y^n_s)^{q+2}} ds)_{n \geq 1}$ converges a.s. and in $L^1$ and we denote by $\Phi$ the limit

$$\Phi_t = \lim_{n \to +\infty} \frac{q(q+1)}{2} \mathbb{E}^{\mathcal{F}_t} \int_{t\wedge\tau}^{\tau} \frac{\|Z^n_s\|^2}{(Y^n_s)^{q+2}} ds.$$

On the set $\{t \geq \tau\}$, $\Phi_t = 0$ a.s. and we have

$$\frac{q(q+1)}{2} \mathbb{E}^{\mathcal{F}_t} \int_{t\wedge\tau}^{\tau} \frac{\|Z^n_s\|^2}{(Y^n_s)^{q+2}} ds$$

$$\leq q\mathbb{E}^{\mathcal{F}_t}(\tau - t \wedge \tau) + \mathbb{E}^{\mathcal{F}_t}\left(\frac{1}{(\xi \wedge n)^q}\right)$$

$$\leq q\mathbb{E}^{\mathcal{F}_t}(\tau) + \frac{1}{\alpha^q}.$$

Thus,

$$\Phi_t \leq q\mathbb{E}^{\mathcal{F}_t}(\tau) + \frac{1}{\alpha^q}.$$

For $r \leq t$,

$$\int_{r \wedge \tau}^{\tau} \frac{\|Z_s^n\|^2}{(Y_s^n)^{q+2}} ds \geq \int_{t \wedge \tau}^{\tau} \frac{\|Z_s^n\|^2}{(Y_s^n)^{q+2}} ds,$$

$$\implies \mathbb{E}^{\mathcal{F}_r} \int_{r \wedge \tau}^{\tau} \frac{\|Z_s^n\|^2}{(Y_s^n)^{q+2}} ds \geq \mathbb{E}^{\mathcal{F}_r} \mathbb{E}^{\mathcal{F}_t} \int_{t \wedge \tau}^{\tau} \frac{\|Z_s^n\|^2}{(Y_s^n)^{q+2}} ds,$$

$$\implies \Phi_r \geq \mathbb{E}^{\mathcal{F}_r} \Phi_t.$$

We deduce that $(\Phi_t)_{0 \leq t}$ is a nonnegative supermartingale. Now for all $n \in \mathbb{N}^*$,

$$\frac{1}{(Y_t^n)^q} = q\mathbb{E}^{\mathcal{F}_t}(\tau - t \wedge \tau) + \mathbb{E}^{\mathcal{F}_t}\left(\frac{1}{(\xi \wedge n)^q}\right) - \frac{q(q+1)}{2}\mathbb{E}^{\mathcal{F}_t}\int_{t \wedge \tau}^{\tau} \frac{\|Z_s^n\|^2}{(Y_s^n)^{q+2}} ds.$$

Fix $t \geq 0$. Taking the limit as $n \to +\infty$, we deduce

$$\frac{1}{(Y_{t \wedge \tau})^q} = q\mathbb{E}^{\mathcal{F}_t}(\tau - t \wedge \tau) + \mathbb{E}^{\mathcal{F}_t}\left(\frac{1}{\xi^q}\right) - \Phi_t.$$

This achieves the proof of Proposition 6. $\square$

$\Phi$ being a nonnegative supermartingale, the limit of $\Phi_{t \wedge \tau}$ as $t$ goes to $+\infty$ exists **P**-a.s. and this limit $\Phi_{\tau^-}$ is finite **P**-a.s. The $L^1$-bounded martingale $\mathbb{E}^{\mathcal{F}_t}(\frac{1}{\xi^q})$ converges a.s. to $1/\xi^q$ as $t$ goes to $+\infty$, then the limit of $Y_{t \wedge \tau}$ as $t \to +\infty$ exists and is equal to

$$\lim_{t \to +\infty} Y_{t \wedge \tau} = \frac{1}{(1/\xi^q - \Phi_{\tau^-})^{1/q}}.$$

If we were able to prove that $\Phi$ is continuous (or $\Phi_{\tau^-}$ is zero a.s.), we would have shown that $Y$ is a continuous process.

3.1.2. *The case $\xi$ nonnegative.* Now we just assume that $\xi \geq 0$. We cannot apply the same arguments because $Y^n$ may to equal to zero with positive probability, which implies, in particular, that $(Y_t^n)^{-q} \notin L^1(\Omega)$. We will approach $Y^n$ in the following way. We define for $n \geq 1$ and $m \geq 1$, $\xi^{n,m}$ by

$$\xi^{n,m} = (\xi \wedge n) \vee \frac{1}{m}.$$

This random variable is in $L^2$ and is greater or equal to $1/m$ a.s. The BSDE (3) with $\xi^{n,m}$ as terminal condition has a unique solution $(\widetilde{Y}^{n,m}, \widetilde{Z}^{n,m})$. It is immediate that if $m \leq m'$ and $n \leq n'$, then

$$\widetilde{Y}^{n,m'} \leq \widetilde{Y}^{n',m}.$$

As for the sequence $Y^n$, we can define $\widetilde{Y}^m$ as the limit when $n$ grows to $+\infty$ of $Y^{n,m}$. That limit $\widetilde{Y}^m$ is greater than $Y = \lim_{n \to +\infty} Y^n$. But for $m \leq m'$ for $0 \leq t \leq T$,

$$\widetilde{Y}^{n,m}_{t\wedge\tau} - \widetilde{Y}^{n,m'}_{t\wedge\tau} = \widetilde{Y}^{n,m}_{T\wedge\tau} - \widetilde{Y}^{n,m'}_{T\wedge\tau} - \int_{t\wedge\tau}^{T\wedge\tau} [(\widetilde{Y}^{n,m}_r)^{q+1} - (\widetilde{Y}^{n,m'}_r)^{q+1}] dr$$

$$- \int_{t\wedge\tau}^{T\wedge\tau} [\widetilde{Z}^{n,m}_r - \widetilde{Z}^{n,m'}_r] dB_r$$

$$\leq \widetilde{Y}^{n,m}_{T\wedge\tau} - \widetilde{Y}^{n,m'}_{T\wedge\tau} - \int_{t\wedge\tau}^{T\wedge\tau} [\widetilde{Z}^{n,m}_r - \widetilde{Z}^{n,m'}_r] dB_r$$

and taking the conditional expectation given $\mathcal{F}_t$,

$$0 \leq \widetilde{Y}^{n,m}_{t\wedge\tau} - \widetilde{Y}^{n,m'}_{t\wedge\tau} \leq \mathbb{E}^{\mathcal{F}_t}(\widetilde{Y}^{n,m}_{T\wedge\tau} - \widetilde{Y}^{n,m'}_{T\wedge\tau}) \leq \frac{1}{m}.$$

Letting first $T \to +\infty$ and then $m' \to +\infty$ in the last estimate leads to

$$0 \leq \widetilde{Y}^{n,m}_{t\wedge\tau} - Y^n_{t\wedge\tau} \leq \frac{1}{m}.$$

Therefore, **P**-a.s.,

$$\sup_{t\geq 0} |\widetilde{Y}^m_{t\wedge\tau} - Y_{t\wedge\tau}| \leq \frac{1}{m}.$$

Since for each $m \geq 0$, $(\widetilde{Y}^m_{t\wedge\tau})_{t\geq 0}$ has a limit on the left at $+\infty$, so does $Y$.

3.2. *Continuity of $Y$.* We know now that

(24) $$\liminf_{t\to+\infty} Y_{t\wedge\tau} = \lim_{t\to+\infty} Y_{t\wedge\tau} \geq \xi$$

and on the set $\{t \geq \tau\}$, $Y_t = \xi$. In this part we give sufficient conditions to ensure that the process $Y$ is continuous, that is,

$$\lim_{t\to+\infty} Y_{t\wedge\tau} = \xi.$$

It suffices to prove the result on the set $\{\xi < \infty\}$. In the rest of this section, we will suppose that (A1) and (A2) hold, and $\mathbf{P}(\xi < \infty) > 0 \Rightarrow F_\infty \neq \partial D$.

3.2.1. *A first step.* In the first section we have proved the following estimate:

(12) $\quad$ **P**-a.s. $\quad \forall t \geq 0 \ |Y_t| \leq \dfrac{C}{(\rho(X^x_{t\wedge \tau_x}))^{2/q}},$

where $\rho$ is the distance to the boundary of $D$. The constant $C$ depends on $q$, $D$ and the bound on $b$ and $\sigma$ in (B). Here we want to construct another estimate which depends also on the function $g$. Our result is the following:

PROPOSITION 7. *Suppose that the boundary of $D$ belongs to $C^3$. If $U$ is an open set such that $\overline{U} \cap F_\infty = \varnothing$ and $U \cap \partial D \neq \varnothing$, then there exists a constant $C = C(U, g, q, b, \sigma, D)$ and an open set $D_U$ such that $D \subset D_U$ and if $\rho_U$ denotes the distance to the boundary of $D_U$, we have*

(25) $\quad$ **P**-*a.s.* $\quad \forall n \in \mathbb{N}, \ \forall t \geq 0 \ Y^n_t \leq \dfrac{C}{(\rho_U(X_{t\wedge \tau}))^{2/q}}.$

*Recall that $\tau$ is always the first exit time from $\overline{D}$.*

PROOF. We suppose that the set $F_\infty = \{g = +\infty\}$ is not equal to $\partial D$. Hence, if we define for all $\varepsilon \geq 0$ the set

$$F_\varepsilon = \{y \in \partial D; \ \mathrm{dist}(y, F_\infty) \leq \varepsilon\},$$

there exists $\varepsilon' > 0$ such that $F_{\varepsilon'} \neq \partial D$. Moreover, if $U$ is an open subset of $\mathbb{R}^d$ such that $F_\infty \cap \overline{U} = \varnothing$ and $U \cap \partial D \neq \varnothing$, there exists $0 < \varepsilon < \varepsilon'$ such that $U \cap \partial D \subset \partial D \setminus F_\varepsilon$.

Recall that $D$ is a bounded open set of $\mathbb{R}^d$ with a boundary $\partial D \in C^3$. Thus, there exists $r > 0$ such that on $\Gamma_r = \{y \in \mathbb{R}^d; \ \mathrm{dist}(y, \partial D) < r\}$, the signed distance $d$

$$d(x) = \begin{cases} \mathrm{dist}(x, \partial D) = \rho(x), & \text{if } x \in D, \\ -\mathrm{dist}(x, \partial D), & \text{if } x \in \mathbb{R}^d \setminus D, \end{cases}$$

belongs to $C^3(\Gamma_r)$. Moreover, for all $y \in \Gamma_r$, there exists a unique $x \in \partial D$ such that $y = x - d(y)\overrightarrow{n}(x)$, where $\overrightarrow{n}(x)$ is the outward normal vector at the point $x \in \partial D$. We have $\|y - x\| = |d(y)| = \mathrm{dist}(y, \partial D)$. The result can be found in [9].

We take a function $\psi_\varepsilon : \mathbb{R}^d \to [0, 1]$ such that $\psi_\varepsilon$ is of class $C^2(\mathbb{R}^d)$ and $\psi_\varepsilon = 0$ on $F_\varepsilon$ and $\psi_\varepsilon = 1$ on $\partial D \setminus F_{2\varepsilon}$. With this function we define the set $D_\varepsilon$ as follows:

$$D_\varepsilon = D \cup \{y \in \mathbb{R}^d; \exists x \in \partial D, \exists \nu \in [0, r/2[ \text{ s.t. } y = x + \nu \psi_\varepsilon(x) \overrightarrow{n}(x)\},$$

where $\overrightarrow{n}(x)$ always denotes the outward normal vector at the point $x \in \partial D$. We can easily prove that $D_\varepsilon$ is included in $D \cup \Gamma_r$, and that

$$F_\varepsilon \subset \partial D_\varepsilon \quad \text{and} \quad \partial D \setminus F_{2\varepsilon} \subset D_\varepsilon.$$

If $\partial D \in C^3$, then the boundary of $D_\varepsilon$ is of class $C^2$, and from our construction, the distance to the boundary of $D_\varepsilon$, denoted by $\rho_\varepsilon$, is also a $C^2$ function on the set $\Gamma_r$. Moreover, if $y \in \partial D \setminus F_{2\varepsilon}$, then $\rho_\varepsilon(y) = r/2 > 0$.

Now the proof of (25) is similar to the proof of Theorem 6. Let $\Phi$ be a $C^\infty(\mathbb{R}^d)$ function such that $\Phi = 1$ on $D \setminus \Gamma_r$ and $\Phi = 0$ on $\Gamma_{r/2}$. For all $0 < \eta < r/2$ and $C > 0$, we define a function $\Psi = \Psi_\eta \in C^2(\mathbb{R}^d; \mathbb{R}_+)$ such that on $\overline{D}$,

$$\Psi = \Psi_\eta = \frac{C}{[(1-\Phi)\rho_\varepsilon + R_\varepsilon \Phi + \eta]^{2/q}} = \frac{C}{[\theta_\eta]^{2/q}}.$$

The constant $R_\varepsilon$ is the supremum of $\rho_\varepsilon$ on $\overline{D}_\varepsilon$ and we can easily see that $R_\varepsilon \leq \sup \rho + r/2$ and that $\theta_\eta \geq \rho_\varepsilon$ for all $\eta > 0$. Remark also that $\Psi$ is of class $C^2$ on $\overline{D}$. We apply the Itô formula to $\Psi(X_{t \wedge \tau})$ and by the same arguments as in the proof of Theorem 6, we can choose the constant $C$ [depending only on $D$, on $q$ and on the bound of $b$ and $\sigma$ in (B)], such that for all $0 \leq t \leq T$,

$$\Psi(X_{t \wedge \tau}) = \Psi(X_{T \wedge \tau}) - \int_{t \wedge \tau}^{T \wedge \tau} \nabla \Psi(X_r) \sigma(X_r) \, dB_r$$
$$- \int_{t \wedge \tau}^{T \wedge \tau} \Psi(X_r)^{1+q} \, dr + \int_{t \wedge \tau}^{T \wedge \tau} U_r \, dr;$$

with $U$ a nonnegative adapted process. The constant $C$ must satisfy (15), that is,

$$C^q + \frac{2\theta_\eta}{q}(\nabla \theta_\eta) b - \frac{1}{q}\left(\frac{2}{q} + 1\right) \|\sigma \nabla \theta_\eta\|^2 + \frac{\theta_\eta}{q} \operatorname{Trace}(\sigma \sigma^* D^2 \theta_\eta) \geq 0.$$

Moreover, on $\{t \geq \tau\}$,

$$\Psi(X_\tau) = \frac{C}{\eta^{2/q}} \quad \text{if } X_\tau \in F_\varepsilon,$$

because on $\partial D$, $\Phi = 0$, and on $F_\varepsilon$, $\rho_\varepsilon = 0$;

$$\Psi(X_\tau) \geq \frac{C}{(\eta + r/2)^{2/q}} \quad \text{if } X_\tau \in \partial D \setminus F_\varepsilon,$$

because on $\partial D$, $0 \leq \rho_\varepsilon \leq r/2$. Recall that for all $n \in \mathbb{N}$, $(Y^n, Z^n)$ is the solution of the BSDE (3) with terminal time $\tau$ and terminal data $g \wedge n$. On the compact set $\overline{\partial D \setminus F_\varepsilon}$, by (A2), the function $g$ is bounded by a constant $K = K_\varepsilon$. We choose $C > 0$ and $0 < \eta < r/2$ such that

$$\frac{C}{\eta^{2/q}} \geq n \quad \text{and} \quad \frac{C}{(\eta + r/2)^{2/q}} \geq K.$$

We can take $C > Kr^{2/q}$ satisfying (15), and $\eta < r/2 \wedge C^{q/2}/n^{q/2}$. Note that $C$ does not depend on $\eta$. Therefore, if we define, for $t \geq 0$,

$$\widetilde{Y}_t = \Psi(X_{t \wedge \tau}) \quad \text{and} \quad \widetilde{Z}_t = \nabla \Psi(X_t) \sigma(X_t) \mathbf{1}_{t < \tau},$$

the process $(\widetilde{Y}, \widetilde{Z})$ satisfies **P**-a.s., for all $0 \leq t \leq T$,

$$\widetilde{Y}_{t\wedge\tau} = \widetilde{Y}_{T\wedge\tau} - \int_{t\wedge\tau}^{T\wedge\tau} \widetilde{Y}_r^{1+q}\,dr + \int_{t\wedge\tau}^{T\wedge\tau} U_r\,dr - \int_{t\wedge\tau}^{T\wedge\tau} \widetilde{Z}_r\,dB_r,$$

$U$ being a nonnegative process, and on the set $\{t \geq \tau\}$: $\widetilde{Y}_t \geq g(X_\tau) \wedge n$. From the comparison theorem (Corollary 4.4.2 in [6]) for solutions of a BSDE, we obtain

**P**-a.s. $\quad \forall t \geq 0 \; Y_t^n \leq \Psi_\eta(X_{t\wedge\tau}) \leq \dfrac{C}{(\rho_\varepsilon(X_{t\wedge\tau}))^{2/q}}.$

Since this inequality holds for all $n$, we have proved the proposition. □

The main interest of Proposition 7 is that if $h \in C_0(U)$ ($h$ has a compact support included in $U$) with $\overline{U} \cap F_\infty = \varnothing$ and $U \cap \partial D \neq \varnothing$, then the sequence $h(X_{\cdot \wedge \tau}) Y^n$ is bounded in $L^\infty([0, +\infty[\times\Omega)$: there exists a constant $K = K_U$ such that

**P**-a.s. $\quad \forall t \geq 0 \; h(X_{t\wedge\tau}) Y_t^n \leq K.$

We can also prove the following:

PROPOSITION 8. *For all $\nu > 1$, there exists a constant $K = K_{U,\nu} > 0$ such that*

$$\mathbb{E}\int_0^\tau \|Z_t\|^2 \rho_U^{4/q+\nu}(X_t)\,dt \leq K.$$

PROOF. Using Proposition 7, the proof is the same as the proof of Proposition 4. □

3.2.2. *Continuity: the conclusion.* Recall that $F_\infty = \{g = +\infty\} \cap \partial D$ is a closed set, that $U$ is an bounded open set such that $\overline{U} \cap F_\infty = \varnothing$ and $U \cap \partial D \neq \varnothing$.

Now we take a function $\varphi : \mathbb{R}^d \to \mathbb{R}_+$ of class $C^2$ and with a compact support included in $U$. For $\beta > 0$, we apply the Itô formula to the process $e^{-\beta t} Y_t^n \varphi(X_t)$:

$$\begin{aligned}
&\mathbb{E}(e^{-\beta\tau}(g \wedge n)(X_\tau)\varphi(X_\tau)) \\
&\quad = \mathbb{E}(e^{-\beta(t\wedge\tau)} Y_{t\wedge\tau}^n \varphi(X_{t\wedge\tau})) \\
&\qquad - \beta \mathbb{E}\int_{t\wedge\tau}^\tau e^{-\beta r}\varphi(X_r) Y_r^n\,dr + \mathbb{E}\int_{t\wedge\tau}^\tau e^{-\beta r}\varphi(X_r)Y_r^n|Y_r^n|^q\,dr \\
&\qquad + \mathbb{E}\int_{t\wedge\tau}^\tau e^{-\beta r} Y_r^n \mathcal{L}\varphi(X_r)\,dr + \mathbb{E}\int_{t\wedge\tau}^\tau e^{-\beta r} Z_r^n \nabla\varphi(X_r)\sigma(X_r)\,dr,
\end{aligned}$$

(26)

where $\mathcal{L}$ is defined by (7). In (26), every term, except maybe the last one, is well defined because $\beta > 0$, $\varphi$ is a $C^2$ function with compact support, and $Y^n$ is bounded by $n$. Now using the Cauchy–Schwarz inequality, we obtain, for all $\eta > 1$,

$$\mathbb{E} \int_0^\tau e^{-\beta r} |Z_r^n \cdot \nabla\varphi(X_r)\sigma(X_r)| \, dr$$
$$(27) \qquad \leq \left[ \mathbb{E} \int_0^\tau \|Z_r^n\|^2 \rho_U^{4/q+\eta}(X_r) \, dr \right]^{1/2}$$
$$\times \left[ \mathbb{E} \int_0^\tau e^{-2\beta r} \rho_U^{-4/q-\eta}(X_r) \|\nabla\varphi(X_r)\sigma(X_r)\|^2 \, dr \right]^{1/2}.$$

We already know that $\sigma$ is bounded on $\overline{D}$. The support of $\nabla\varphi$ is included in $U$. In our previous construction of $D_U$, we have $U \cap \overline{D} \subset D_U$ and on $U \cap \overline{D}$, $\rho_U \geq r/2 > 0$. Therefore, $\rho_U^{-4/q-\eta} \nabla\varphi$ is a continuous and bounded function. With Proposition 8, we deduce

$$(28) \qquad \mathbb{E} \int_0^\tau e^{-\beta r} |Z_r^n \cdot \nabla\varphi(X_r)\sigma(X_r)| \, dr \leq K.$$

It is important to remark that the constant $K$ does not depend on $n$.

We want to pass to the limit when $n \to +\infty$ in (26). With the monotone convergence theorem, we obtain, for all $0 \leq t$,

$$\lim_{n \to +\infty} \mathbb{E}(e^{-\beta\tau}(g \wedge n)(X_\tau)\varphi(X_\tau)) = \mathbb{E}(e^{-\beta\tau} g(X_\tau)\varphi(X_\tau));$$

$$\lim_{n \to +\infty} \mathbb{E} \int_{t \wedge \tau}^\tau e^{-\beta r} \varphi(X_r) Y_r^n \, dr = \mathbb{E} \int_{t \wedge \tau}^\tau e^{-\beta r} \varphi(X_r) Y_r \, dr;$$

$$\lim_{n \to +\infty} \mathbb{E} \int_{t \wedge \tau}^\tau e^{-\beta r} \varphi(X_r)(Y_r^n)^{1+q} \, dr = \mathbb{E} \int_{t \wedge \tau}^\tau e^{-\beta r} \varphi(X_r)(Y_r)^{1+q} \, dr.$$

The support of the function $\mathcal{L}\varphi$ is included in $U$. Therefore, from Proposition 7, $Y^n \mathcal{L}\varphi(X)$ is a.s. bounded. From the dominated convergence theorem, we deduce that, for all $t \geq 0$,

$$\lim_{n \to +\infty} \mathbb{E} \int_{t \wedge \tau}^\tau e^{-\beta r} Y_r^n \mathcal{L}\varphi(X_r) \, dr = \mathbb{E} \int_{t \wedge \tau}^\tau e^{-\beta r} Y_r \mathcal{L}\varphi(X_r) \, dr.$$

The last term in (26) is equal to

$$\mathbb{E} \int_{t \wedge \tau}^\tau e^{-\beta r} Z_r^n \cdot \nabla\varphi(X_r)\sigma(X_r) \, dr$$
$$= \mathbb{E} \int_{t \wedge \tau}^\tau e^{-\beta r} \rho_U^{2/q+\eta/2}(X_r) Z_r^n \cdot \rho_U^{-2/q-\eta/2}(X_r) \nabla\varphi(X_r)\sigma(X_r) \, dr.$$

From Proposition 8, the sequence $\rho_U^{2/q+\eta/2}(X)Z^n \mathbf{1}_{\tau >.}$ is bounded in $L^2([0, +\infty[\times\Omega)$ for all $\eta > 1$. Therefore, after extraction of a suitable subsequence, which we omit as an abuse of notation, $\rho_U^{2/q+\eta/2}(X)Z^n \mathbf{1}_{\tau>\cdot}$ converges weakly in $L^2([0,+\infty[\times\Omega)$. The process

$$e^{-\beta\cdot}\rho_U^{-2/q-\eta/2}(X)\nabla\varphi(X)\sigma(X)\mathbf{1}_{\tau>\cdot}$$

is in $L^2([0,+\infty[\times\Omega)$ because $\sigma$ is bounded, and from our construction, $\rho_U^{-2/q-\eta/2}(X)\nabla\varphi(X)\mathbf{1}_{\tau>\cdot}$ is also bounded. Using Proposition 3, we obtain

$$\lim_{n\to+\infty} \mathbb{E}\int_{t\wedge\tau}^\tau e^{-\beta r} Z_r^n \cdot \nabla\varphi(X_r)\sigma(X_r)\,dr = \mathbb{E}\int_{t\wedge\tau}^\tau e^{-\beta r} Z_r \cdot \nabla\varphi(X_r)\sigma(X_r)\,dr.$$

Finally, letting $n \to +\infty$ in (26), we have, for all $0 \leq t$,

$$\mathbb{E}(e^{-\beta\tau}g(X_\tau)\varphi(X_\tau))$$
$$= \mathbb{E}(e^{-\beta(t\wedge\tau)}Y_{t\wedge\tau}\varphi(X_{t\wedge\tau}))$$
(29)
$$-\beta\mathbb{E}\int_{t\wedge\tau}^\tau e^{-\beta r}\varphi(X_r)Y_r\,dr + \mathbb{E}\int_{t\wedge\tau}^\tau e^{-\beta r}\varphi(X_r)(Y_r)^{1+q}\,dr$$
$$+ \mathbb{E}\int_{t\wedge\tau}^\tau e^{-\beta r}Y_r\mathcal{L}\varphi(X_r)\,dr + \mathbb{E}\int_{t\wedge\tau}^\tau e^{-\beta r}Z_r \cdot \nabla\varphi(X_r)\sigma(X_r)\,dr.$$

From Proposition 7, we know that

$$\mathbb{E}\int_0^\tau e^{-\beta r}\varphi(X_r)Y_r\,dr \leq K,$$

$$\mathbb{E}\int_0^\tau e^{-\beta r}\varphi(X_r)(Y_r)^{1+q}\,dr + \mathbb{E}\int_0^\tau e^{-\beta r}Y_r|\mathcal{L}\varphi(X_r)|\,dr \leq K.$$

For the last term, using the Cauchy–Schwarz inequality [see (27)] and Proposition 4, we obtain

$$\mathbb{E}\int_0^\tau e^{-\beta r}|Z_r \cdot \nabla\varphi(X_r)\sigma(X_r)|\,dr \leq K.$$

Therefore, when $t$ goes to $+\infty$ in the equation (29), we obtain, using Fatou's lemma,

(30)
$$\mathbb{E}(e^{-\beta\tau}g(X_\tau)\varphi(X_\tau)) = \lim_{t\to+\infty}\mathbb{E}(e^{-\beta(t\wedge\tau)}Y_{t\wedge\tau}\varphi(X_{t\wedge\tau}))$$
$$\geq \mathbb{E}\left[e^{-\beta\tau}\varphi(X_\tau)\left(\lim_{t\to+\infty}Y_{t\wedge\tau}\right)\right].$$

But recall that we already know that

(24)
$$\lim_{t\to+\infty}Y_{t\wedge\tau} \geq g(X_\tau).$$

Hence, the inequality in (30) is in fact an equality, that is,

$$\mathbb{E}(e^{-\beta\tau}g(X_\tau)\varphi(X_\tau)) = \mathbb{E}\left[e^{-\beta\tau}\varphi(X_\tau)\left(\lim_{t\to+\infty}Y_{t\wedge\tau}\right)\right].$$

And using again (24), we conclude that
$$\lim_{t \to +\infty} Y_{t \wedge \tau} = g(X_\tau), \qquad \textbf{P}\text{-a.s. on } \{g(X_\tau) < \infty\}.$$

**4. Minimal solution.** In the third section we have constructed a process $(Y, Z)$ which satisfies the conditions (D1) and (D2) of Definition 2. We will prove now that, if this process is a solution of the BSDE (3), that is, if it satisfies also the condition (D3), then it is the minimal nonnegative solution.

THEOREM 7. *Let the conditions* (L), (B) *and* (E) *hold and let* $(\overline{Y}, \overline{Z})$ *be a nonnegative solution of the BSDE* (3) (*solution in the sense of Definition* 2). *Then* **P**-*a.s. for all* $t \geq 0$,
$$\overline{Y}_t \geq Y_t.$$

PROOF. Recall that $\tau_\eta$ is the first exit time of $D \setminus \Gamma_\eta$ and we have, for all $0 \leq t \leq T$,
$$\overline{Y}_{t \wedge \tau_\eta} = \overline{Y}_{T \wedge \tau_\eta} - \int_{t \wedge \tau_\eta}^{T \wedge \tau_\eta} (\overline{Y}_r)^{1+q} \, dr - \int_{t \wedge \tau_\eta}^{T \wedge \tau_\eta} \overline{Z}_r \, dB_r.$$
For $n \in \mathbb{N}^*$, $(Y^n, Z^n)$ is the solution (in the sense of Definition 1) of the BSDE (3) with terminal data $\xi \wedge n$. We compare $\overline{Y}$ with $Y^n$:
$$\overline{Y}_{t \wedge \tau_\eta} - Y^n_{t \wedge \tau_\eta} = \overline{Y}_{T \wedge \tau_\eta} - Y^n_{T \wedge \tau_\eta} \int_{t \wedge \tau_\eta}^{T \wedge \tau_\eta} (\overline{Y}_r)^{1+q} - (Y^n_r)^{1+q} \, dr$$
$$- \int_{t \wedge \tau_\eta}^{T \wedge \tau_\eta} (\overline{Z}_r - Z^n_r) \, dB_r$$
$$= \overline{Y}_{T \wedge \tau_\eta} - Y^n_{T \wedge \tau_\eta} - \int_{t \wedge \tau_\eta}^{T \wedge \tau_\eta} \alpha^n_r (\overline{Y}_r - Y^n_r) \, dr$$
$$- \int_{t \wedge \tau_\eta}^{T \wedge \tau_\eta} (\overline{Z}_r - Z^n_r) \, dB_r,$$
where the process $\alpha^n_r$ is defined by
$$\alpha^n_r = \frac{(\overline{Y}_r)^{1+q} - (Y^n_r)^{1+q}}{\overline{Y}_r - Y^n_r}, \qquad \text{if } \overline{Y}_r \neq Y^n_r,$$
$$\alpha^n_r = (1+q)(Y^n_r)^q, \qquad \text{if } \overline{Y}_r = Y^n_r.$$
$\alpha^n$ is a nonnegative process and we have a linear BSDE whose solution is
$$(31) \quad \overline{Y}_{t \wedge \tau_\eta} - Y^n_{t \wedge \tau_\eta} = \mathbb{E}^{\mathcal{F}_t}\left[(\overline{Y}_{T \wedge \tau_\eta} - Y^n_{T \wedge \tau_\eta}) \exp\left(-\int_{t \wedge \tau_\eta}^{T \wedge \tau_\eta} \alpha^n_r \, dr\right)\right].$$

From the hypothesis of the theorem, $\overline{Y}$ is nonnegative and $Y^n$ is bounded by $n$. Indeed, on the set $\{t \geq \tau\}$, $Y_t^n = \xi \wedge n \leq n$ and for all $0 \leq t \leq T$,

$$Y_{t\wedge\tau}^n = Y_{T\wedge\tau}^n - \int_{t\wedge\tau}^{T\wedge\tau} (Y_r^n)^{1+q}\, dr - \int_{t\wedge\tau}^{T\wedge\tau} Z_r^n\, dB_r$$

$$\leq Y_{T\wedge\tau}^n - \int_{t\wedge\tau}^{T\wedge\tau} Z_r^n\, dB_r,$$

thus,

$$Y_{t\wedge\tau}^n \leq Y_\tau^n - \int_{t\wedge\tau}^\tau Z_r^n\, dB_r = \xi \wedge n - \int_{t\wedge\tau}^\tau Z_r^n\, dB_r \leq n - \int_{t\wedge\tau}^\tau Z_r^n\, dB_r.$$

Taking the conditional expectation, we deduce that, for all $t \geq 0$, $Y_t^n \leq n$.

We now pass to the limit in (31) first as $\eta \to 0$, then as $T \to +\infty$ and with the Fatou lemma, we obtain for all $t \geq 0$,

$$\overline{Y}_{t\wedge\tau} - Y_{t\wedge\tau}^n \geq 0.$$

Therefore, $\overline{Y}$ is greater than $Y^n$ for all $n \in \mathbb{N}^*$ and thus greater than $Y$. □

Moreover, we obtain the following result:

PROPOSITION 9. *There exists a constant $C$ (the same constant as in Theorem* 6*) such that $\mathbf{P}$-a.s., for all $t \geq 0$,*

$$\overline{Y}_t \leq \frac{C}{\rho^{2/q}(X_{t\wedge\tau})}.$$

PROOF. For all sufficiently small $\eta > 0$, we denote by $\rho_\eta$ the distance from the boundary of $D_\eta = D \setminus \Gamma_\eta$, that is,

$$D_\eta = \{x \in D, \rho(x) \geq \eta\}.$$

If $x \in D_\eta$, $\rho(x) - \eta \leq \rho_\eta(x) \leq \rho(x)$. We consider the first exit time

$$\tau_\eta = \inf\{t \geq 0, X_t \notin \overline{D}_\eta\}.$$

From Theorem 6 we deduce

$$\forall t \geq 0 \qquad \overline{Y}_{t\wedge\tau_\eta} \leq \frac{C}{\rho_\eta^{2/q}(X_{t\wedge\tau_\eta})} \leq \frac{C}{\rho^{2/q}(X_{t\wedge\tau_\eta}) - \eta}.$$

The constant $C$ which appears in the previous inequality may depend on $\eta$. In the proof of Theorem 6 we use the fact that there exists $\mu > 0$ such that on $\Gamma_\mu$, the signed distance function is of class $C^2$. But if $\eta < \mu$, it is also true that $\rho_\eta$ is of class $C^2$ on $\Gamma_\mu$. So in the proof of the theorem we can use the same function $\varphi$ and the same bound $R$ for $\rho_\eta$ and $\rho$. Moreover, on $\Gamma_\mu$, $|\nabla \rho| = 1$ and $D^2 \rho$ depends only on the curvature of $\partial D$. Therefore, we can choose a constant $C$ independent of $\eta$ if $\eta < \mu$.

To conclude, let $\eta \to 0$ and we obtain the desired inequality. □

**5. Viscosity solution of the associated elliptic PDE.** Recall that $D$ is a bounded open subset of $\mathbb{R}^d$ with a $C^3$ boundary. For all $x \in \overline{D}$, $\{X^x_t; t \geq 0\}$ is the solution of the SDE (4):

$$(4) \qquad X^x_t = x + \int_0^t b(X^x_r)\, dr + \int_0^t \sigma(X^x_r)\, dB_r \qquad \text{for } t \geq 0.$$

The functions $b$ and $\sigma$ are defined on $\mathbb{R}^d$, with values respectively in $\mathbb{R}^d$ and $\mathbb{R}^{d \times d}$, and such that $b$ and $\sigma$ are continuous on $\mathbb{R}^d$ and satisfy the conditions (M), (L) and (B). For each $x \in \overline{D}$, we define the stopping time $\tau_x = \inf\{t \geq 0, X^x_t \notin \overline{D}\}$. We assume that

$$(32) \qquad \mathbf{P}(\tau_x < \infty) = 1 \qquad \text{for all } x \in \overline{D},$$

that the set of singular points

$$(\text{C1}) \qquad \Gamma = \{x \in \partial D; \mathbf{P}(\tau_x > 0) > 0\} \qquad \text{is empty},$$

and that for some $\beta > 0$ and all $x \in \overline{D}$,

$$(\text{C2}) \qquad \mathbb{E} e^{\beta \tau_x} < \infty.$$

Let us recall the following result (cf. Proposition 5.2. in [18]):

PROPOSITION 10. *Under the conditions* (C1) *and* (C2), *the mapping* $x \mapsto \tau_x$ *is a.s. continuous on* $\overline{D}$.

Let $g: \partial D \to \overline{\mathbb{R}}_+$ be a continuous function and for all $n \in \mathbb{N}$, we define $g_n = g \wedge n$. Hence, $g_n$ is a continuous function. For all $n \in \mathbb{N}$, from Remark 4, $\{(Y^{x,n}_t, Z^{x,n}_t); t \geq 0\}$ is the unique solution (in the sense of Definition 1) of the BSDE (3)

$$(33) \qquad Y^{x,n}_t = g_n(X^x_{\tau_x}) - \int_{t \wedge \tau_x}^{\tau_x} Y^{x,n}_r |Y^{x,n}_r|^q\, dr - \int_{t \wedge \tau_x}^{\tau_x} Z^{x,n}_r\, dB_r.$$

We denote by $u_n$ the function defined on $\overline{D}$ by

$$u_n(x) \triangleq Y^{x,n}_0.$$

For $h \in C(\partial D, \mathbb{R})$, we consider the elliptic PDE (6) with boundary condition $h$:

$$-\mathcal{L}v + v|v|^q = 0 \qquad \text{on } D;$$
$$v = h \qquad \text{on } \partial D.$$

The following definition can be found in [1] and [2] (or [5] and [18] for $v$ continuous). If $v$ is a function defined on $\overline{D}$, we denote by $v^*$ (resp. $v_*$) the upper- (resp. lower-) semicontinuous envelope of $v$: for all $x \in \overline{D}$,

$$v^*(x) = \limsup_{x' \to x, x' \in \overline{D}} v(x') \quad \text{and} \quad v_*(x) = \liminf_{x' \to x,\ x' \in \overline{D}} v(x').$$

DEFINITION 3 (*Viscosity solution*).

- $v : \overline{D} \to \mathbb{R}$ is called a *viscosity subsolution* of (6) if $v^* < +\infty$ on $\overline{D}$ and if for all $\phi \in C^2(\mathbb{R}^d)$, whenever $x \in \overline{D}$ is a point of local maximum of $v^* - \phi$,

$$-\mathcal{L}\phi(x) + v^*(x)|v^*(x)|^q \leq 0 \quad \text{if } x \in D;$$

$$\min(-\mathcal{L}\phi(x) + v^*(x)|v^*(x)|^q, v^*(x) - h(x)) \leq 0 \quad \text{if } x \in \partial D.$$

- $v : \overline{D} \to \mathbb{R}$ is called a *viscosity supersolution* of (6) if $v_* > -\infty$ on $\overline{D}$ and if for all $\phi \in C^2(\mathbb{R}^d)$, whenever $x \in \overline{D}$ is a point of local minimum of $v_* - \phi$,

$$-\mathcal{L}\phi(x) + v_*(x)|v_*(x)|^q \geq 0 \quad \text{if } x \in D;$$

$$\max(-\mathcal{L}\phi(x) + v_*(x)|v_*(x)|^q, v(x) - h(x)) \geq 0 \quad \text{if } x \in \partial D.$$

- $v : \overline{D} \to \mathbb{R}$ is called a *viscosity solution* of (6) if it is both a viscosity sub- and supersolution.

Let us recall the following result (cf. Theorem 5.3. in [18]):

THEOREM 8. *Under the assumptions* (M), (L), (B), (C1) *and* (C2), *since* $g \wedge n$ *is continuous on* $\partial D$, $u_n$ *is continuous on* $\overline{D}$ *and it is a viscosity solution of the elliptic PDE* (6) *with boundary data* $g \wedge n$.

REMARK 6. Since $g \wedge n$ is continuous on $\partial D$, from Theorem 3.3 in [5], it follows that $u_n$ is the unique continuous viscosity solution of the PDE (6) with terminal data $g \wedge n$.

From now on we add the uniformly elliptic condition: there exists a constant $\alpha > 0$ such that, for all $x \in \mathbb{R}^d$,

(E) $$\sigma \sigma^*(x) \geq \alpha \mathrm{Id}.$$

With this assumption, (C1) and (C2) hold if (B) is true. In the previous sections we have constructed a process $\{(Y^x_t, Z^x_t); t \geq 0\}$ which is a solution of the BSDE (3) with terminal data $g(X^x_{\tau_x})$ (in the sense of Definition 2). $Y^x$ is the limit of $Y^{x,n}$: for all $t \geq 0$,

(16) $$Y^x_t = \lim_{n \to +\infty} Y^{x,n}_t.$$

If we define

$$u(x) \triangleq Y^x_0,$$

then $u$ is the limit of the sequence $u_n$. Thus, $u$ is nonnegative. Since $u$ is the supremum of continuous functions $u_n$, $u$ is lower-semicontinuous on $\overline{D}$ and satisfies

(34) $$\forall x \in \overline{D} \quad u(x) \leq \frac{C}{\rho^{2/q}(x)}.$$

Recall that $\rho$ is the distance from the boundary $\partial D$ and $C$ is a constant which does not depend on $g$. Moreover, $u(x) = g(x)$ on $\partial D$. Since $g$ is not bounded on $\partial D$, we cannot apply Theorem 8. Moreover, the condition $v^* < +\infty$ in Definition 3 cannot be satisfied on $\overline{D}$. Therefore, we change the definition of a solution.

DEFINITION 4 (*Unbounded viscosity solution*). We say that $v$ is a viscosity solution of the PDE

$$-\mathcal{L}v + v|v|^q = 0 \quad \text{on } D,$$
$$v = g \quad \text{on } \partial D,$$
(6)

with unbounded terminal data $g$ if $v$ is a viscosity solution on $D$ in the sense of Definition 3 and if

$$g(x) \leq \lim_{\substack{x' \to x \\ x' \in D, x \in \partial D}} v_*(x') \leq \lim_{\substack{x' \to x \\ x' \in D, x \in \partial D}} v^*(x') \leq g(x).$$

Remark that this definition implies that $v^* < +\infty$ and $v_* > -\infty$ on $D$.

### 5.1. *u is a viscosity solution.*

LEMMA 4. *The function $u$ is a viscosity solution of the PDE (6) on $D$.*

PROOF. We will use the half-relaxed upper- and lower-limit of the sequence of functions $u_n$:

$$\bar{u}(x) = \limsup_{\substack{n \to +\infty \\ x' \to x}} u_n(x') \quad \text{and} \quad \underline{u}(x) = \liminf_{\substack{n \to +\infty \\ x' \to x}} u_n(x').$$

Since $\{u_n\}$ is a nondecreasing sequence of continuous functions, we have

$$\forall x \in D \quad u(x) = u_*(x) = \underline{u}(x) \leq u^*(x) = \bar{u}(x).$$

We fix $\eta > 0$ and we prove that on $\widetilde{D} = D \setminus \{x \in D, \rho(x) \leq \eta\}$, $u$ is a viscosity solution. We already know that

$$\forall x \in \widetilde{D} \quad \bar{u}(x) \leq \frac{C}{\eta^{2/q}}.$$

Recall that $u_n$ is a continuous viscosity solution and from the Lemma 6.1 of [5], we deduce that $u$ is a viscosity solution of (6) on $\widetilde{D}$ and this holds for all $\eta > 0$. Therefore, the lemma is proved. □

Since $u_n$ is a nondecreasing sequence of $C^0(\overline{D})$ functions, we have,

(35) $$\forall x \in \partial D \quad \liminf_{x' \to x, \ x' \in D} u(x') \geq g(x) = u(x).$$

Hence, $u_*$ is a supersolution of (6) because $u_* \geq g$ on $\partial D$.

BSDE WITH SINGULAR FINAL CONDITION 39Lemma 5. *The solution $u$ satisfies the boundary condition, that is,*

$$\lim_{\substack{x' \to x \\ x' \in D,\ x \in \partial D}} u^*(x') \leq g(x) = u(x).$$

Proof. We already know that

$$\liminf_{\substack{x' \to x \\ x' \in D, x \in \partial D}} u(x') \geq g(x) = u(x).$$

So we just have to prove the converse inequality on the set $\{g < +\infty\}$. If $U$ is an open set such that $\overline{U} \cap F_\infty = \varnothing$ and $U \cap \partial D \neq \varnothing$, there exists an open set $D_U$ and a constant $C_U$ such that, for all $n \in \mathbb{N}$,

$$\mathbf{P}\text{-a.s.} \qquad \forall t \geq 0\ Y_t^{x,n} \leq \frac{C_U}{\rho_U^{2/q}(X_{t \wedge \tau_x}^x)}.$$

Recall that $\rho_U$ is the distance to the boundary of $D_U$. From the proof of Proposition 7, the choice of the set $D_U$ and of the constant $C_U$ does not depend on $x \in \overline{D}$.

We write again equation (26):

(36)
$$\begin{aligned}
&\mathbb{E}(e^{-\beta \tau_x}(g \wedge n)(X_{\tau_x}^x)\varphi(X_{\tau_x}^x)) \\
&= u_n(x)\varphi(x) \\
&\quad - \beta \mathbb{E}\int_0^{\tau_x} e^{-\beta r}\varphi(X_r^x)Y_r^{x,n}\,dr + \mathbb{E}\int_0^{\tau_x} e^{-\beta r}\varphi(X_r^x)Y_r^{x,n}|Y_r^{x,n}|^q\,dr \\
&\quad + \mathbb{E}\int_0^{\tau_x} e^{-\beta r}Y_r^{x,n}\mathcal{L}\varphi(X_r^x)\,dr + \mathbb{E}\int_0^{\tau_x} e^{-\beta r}Z_r^{x,n}\cdot\nabla\varphi(X_r^x)\sigma(X_r^x)\,dr.
\end{aligned}$$

The function $\varphi : \mathbb{R}^d \to \mathbb{R}_+$ is of class $C^2$ and has a compact support included in $U$. The constant $\beta$ is positive. From Proposition 7, there exists a constant $K_U$ such that, for all $n \in \mathbb{N}$,

$$\left|\mathbb{E}\int_0^{\tau_x} e^{-\beta r}[\varphi(X_r^x)Y_r^{x,n} + \varphi(X_r^x)Y_r^{x,n}|Y_r^{x,n}|^q + Y_r^{x,n}\mathcal{L}\varphi(X_r^x)]\,dr\right|$$
$$\leq K_U \mathbb{E}\int_0^{\tau_x} e^{-\beta r}\,dr.$$

Moreover, using (27), we have

$$\mathbb{E}\int_0^{\tau_x} e^{-\beta r}|Z_r^{x,n}\nabla\varphi(X_r^x)\sigma(X_r^x)|\,dr$$
$$\leq \left[\mathbb{E}\int_0^{\tau_x}\|Z_r^{x,n}\|^2 \rho_U^{4/q+\eta}(X_r^x)\,dr\right]^{1/2}$$
$$\times \left[\mathbb{E}\int_0^{\tau_x} e^{-2\beta r}\rho_U^{-4/q-\eta}(X_r^x)\|\nabla\varphi(X_r^x)\sigma(X_r^x)\|^2\,dr\right]^{1/2}.$$

Recall that $\rho_U^{-4/q-\eta}\nabla\varphi$ is a continuous and bounded function on $\overline{D}$. With Proposition 8, we obtain that

$$\mathbb{E}\int_0^{\tau_x} e^{-\beta r}|Z_r^{x,n}\cdot\nabla\varphi(X_r^x)\sigma(X_r^x)|\,dr \leq K_U\mathbb{E}\int_0^{\tau_x} e^{-\beta r}\,dr.$$

Since $x \mapsto \tau_x$ is a continuous function on $\overline{D}$ and since $\tau_x = 0$ if $x \in \partial D$, we have

$$\lim_{\substack{x'\to x \\ x'\in D, x\in\partial D}} \mathbb{E}\int_0^{\tau_{x'}} e^{-\beta r}\,dr = 0.$$

If $x \in \partial D$ and if $(x_m)_{m\in\mathbb{N}}$ is a sequence of elements of $D$ which converges to $x$, we replace in (36) $n$ by $m$ and $x$ by $x_m$ and we take the limit as $m \to +\infty$. We obtain, by Fatou's lemma,

$$\limsup_{m\to+\infty} u_m(x_m)\varphi(x_m) = \limsup_{m\to+\infty}\mathbb{E}(e^{-\beta\tau_{x_m}}(g\wedge m)(X_{\tau_{x_m}}^{x_m})\varphi(X_{\tau_{x_m}}^{x_m}))$$

$$\leq \mathbb{E}\left(\limsup_{m\to+\infty}[e^{-\beta\tau_{x_m}}(g\wedge m)(X_{\tau_{x_m}}^{x_m})\varphi(X_{\tau_{x_m}}^{x_m})]\right),$$

because $g\varphi$ is a bounded function. By continuity of $x \mapsto X_{\tau_x}^x$ and of $g\varphi$, we have

$$\limsup_{m\to+\infty} u_m(x_m)\varphi(x) = \limsup_{m\to+\infty} u_m(x_m)\varphi(x_m) \leq g(x)\varphi(x).$$

Finally, on $\{g < \infty\}$, we have

$$\limsup_{\substack{x'\to x \\ x'\in D,\ x\in\partial D}} u^*(x') \leq g(x).$$

With inequality (35), this achieves the proof of the lemma. $\square$

5.2. *Some regularity results on u.* We want to prove now that $u$ is continuous on $\overline{D}$. Here it seems to be necessary to assume the condition (E).

First we prove that, under stronger assumptions on $b$ and $\sigma$, $u$ belongs to $C^0(\overline{D};\mathbb{R}_+) \cap C^2(D;\mathbb{R}_+)$.

PROPOSITION 11. *Recall that $b$ and $\sigma$ satisfy always* (L), (B), (E) *and $\partial D \in C^3$. We assume that $b$ and $\sigma$ belong to $C^1(D)$. Then $u$ is in $C^2(D;\mathbb{R}_+)$.*

In order to prove this result, we need the following lemma:

LEMMA 6. *The assumptions of Proposition 11 hold. We consider a continuous function $h: \partial D \to \mathbb{R}$. Suppose that $(Y,Z)$ is the solution (in the sense of Definition 1) of the BSDE*

$$Y_t = h(X_\tau) - \int_{t\wedge\tau}^\tau Y_r|Y_r|^q\,dr - \int_{t\wedge\tau}^\tau Z_r\,dB_r.$$

BSDE WITH SINGULAR FINAL CONDITION 41

*Then there exists a function* $v : \overline{D} \to \mathbb{R}^+$ *of class* $C^0(\overline{D}) \cap C^2(D)$ *such that*

$$\forall t \geq 0 \qquad Y_t = v(X_{t \wedge \tau}) \quad \text{and} \quad Z_t = \nabla v(X_{t \wedge \tau}) \sigma(X_{t \wedge \tau}) \mathbf{1}_{t < \tau}.$$

*Moreover, $v$ is solution of the PDE* (6) *with boundary condition $h$.*

PROOF. Since $h$ is continuous and since $b$ and $\sigma$ belongs to $C^1(D)$, from the Theorem 15.18 in [9], there exists a unique solution $v \in C^0(\overline{D}) \cap C^2(D)$ of the PDE (6) (see also [15]). We prove that, for all $t \geq 0$, $Y_t = v(X_{t \wedge \tau})$ and $Z_t = \nabla v(X_{t \wedge \tau}) \sigma(X_{t \wedge \tau}) \mathbf{1}_{t < \tau}$. We want to apply the Itô formula to the process $v(X)$. But we just have $v \in C^2(D)$ and we do not know if we can define a function $\widetilde{v} \in C^2(\mathbb{R}^d)$ such that $\widetilde{v} = v$ on $D$.

We will use some arguments of the proof of Theorem 15.18 in [9]. We define a sequence $\{h_m\}$ of functions such that $\{h_m\}$ approximates $h$ uniformly on $\partial D$ and $h_m \in C^{2,\gamma}(\overline{D})$. From Theorem 15.10 in [9], there exists a function $v_m$ such that $v_m$ solves the Dirichlet problem (6) with condition $h_m$ on the boundary and $v_m \in C^{2,\gamma}(\overline{D})$. We apply the Itô formula to the process $v_m(X)$: for all $0 \leq t \leq T$,

$$\begin{aligned}
v_m(X_{t \wedge \tau}) &= v_m(X_{T \wedge \tau}) - \int_{t \wedge \tau}^{T \wedge \tau} (\mathcal{L} v_m)(X_r) \, dr \\
&\quad - \int_{t \wedge \tau}^{T \wedge \tau} \nabla v_m(X_r) \sigma(X_r) \, dB_r \\
&= v_m(X_{T \wedge \tau}) - \int_{t \wedge \tau}^{T \wedge \tau} v_m(X_r) |v_m(X_r)|^q \, dr \\
&\quad - \int_{t \wedge \tau}^{T \wedge \tau} \nabla v_m(X_r) \sigma(X_r) \mathbf{1}_{r < \tau} \, dB_r.
\end{aligned} \tag{37}$$

We denote by $(Y^m, Z^m)$ the solution of the BSDE (3) with terminal data $h_m(X_\tau) \in L^\infty(\Omega)$. Uniqueness of solution of this BSDE implies

$$\forall t \geq 0 \qquad Y^m_t = v_m(X_{t \wedge \tau}) \quad \text{and} \quad Z^m_t = \nabla v_m(X_t) \sigma(X_t) \mathbf{1}_{t < \tau}.$$

From (C2) and Remark 4, we know that there exists a constant $C$ such that, for all $m \in \mathbb{N}$,

$$0 \leq Y^m_0 = v_m(x) \leq \mathbb{E}\left[\sup_{t \in [0,\tau]} e^{\beta t} |Y^m_t|^2\right] \leq C \mathbb{E}[e^{\beta \tau} |h_m(X_\tau)|^2].$$

Since $h_m$ converges uniformly to $h$ on $\partial D$, $h_m$ is a bounded sequence in $L^\infty(\partial D)$. Therefore, the sequence $\{v_m\}$ is uniformly bounded on $\overline{D}$.

From Theorems 6.1, 13.1 and 15.3 in [9], the sequence $\{v_m\}$ converges uniformly on compact subsets of $D$, together with its first and second derivatives, to the function $v$.

Since $\{v_m\}$ is uniformly bounded on $\overline{D}$ and converges to $v$, for all $0 \le t \le T$,

(38) $$\lim_{m \to +\infty} Y_t^m = \lim_{m \to +\infty} v_m(X_{t \wedge \tau}) = v(X_{t \wedge \tau}),$$

$$\lim_{m \to +\infty} \int_{t \wedge \tau}^{T \wedge \tau} v_m(X_r)|v_m(X_r)|^q \, dr = \int_{t \wedge \tau}^{T \wedge \tau} v(X_r)|v(X_r)|^q \, dr.$$

Using Itô's formula and the Burkholder–Davis–Gundy inequality, we obtain, for a constant $c$ independent of $m$,

$$\mathbb{E}\left[\sup_{0 \le t \le \tau} |Y_t^m - Y_t|^2 + \int_0^\tau \|Z_t^m - Z_t\|^2 \, dt\right] \le c\mathbb{E}|h_m(X_\tau) - h(X_\tau)|^2.$$

Therefore, with (38), we conclude that a.s. $Y_t = v(X_{t \wedge \tau})$ for all $t \ge 0$. Moreover, there exists a constant $K$ such that

$$\mathbb{E}\int_0^\tau \|\nabla v_m(X_r)\sigma(X_r)\mathbf{1}_{r<\tau}\|^2 \, dr = \mathbb{E}\int_0^\tau \|Z_r^m\|^2 \, dr \le K < +\infty;$$

and with (37) and (38), for all $0 \le t \le T$,

$$\lim_{m \to +\infty} \int_{t \wedge \tau}^{T \wedge \tau} \nabla v_m(X_r)\sigma(X_r)\mathbf{1}_{r \le \tau} \, dB_r = \int_{t \wedge \tau}^{T \wedge \tau} Z_r \, dB_r.$$

Let $\mathcal{K}$ be a compact subset of $D$. Since the first derivatives of $v_m$ converge uniformly on $\mathcal{K}$, from the dominated convergence theorem, we deduce

$$\lim_{m \to +\infty} \int_{t \wedge \tau}^{T \wedge \tau} \|(\nabla v_m(X_r)\sigma(X_r) - \nabla v(X_r)\sigma(X_r))\mathbf{1}_{r<\tau}\mathbf{1}_{\mathcal{K}}(X_r)\|^2 \, dr = 0.$$

Therefore, for all compact subset $\mathcal{K}$ of $D$, **P**-a.s.,

$$Z_t \mathbf{1}_{\mathcal{K}}(X_t) = \nabla v(X_t)\sigma(X_t)\mathbf{1}_{t<\tau}\mathbf{1}_{\mathcal{K}}(X_t).$$

If $\{\mathcal{K}_m\}$ is an increasing sequence of compact subsets of $D$ such that $\bigcup_m \mathcal{K}_m = D$, for all $m$,

$$\mathbb{E}\int_0^\tau \|\nabla v(X_t)\sigma(X_t)\mathbf{1}_{t<\tau}\mathbf{1}_{\mathcal{K}_m}(X_t)\|^2 \, dt \le \mathbb{E}\int_0^\tau \|Z_t\|^2 \, dt < +\infty$$

and since $\tau > t$ implies $X_t \in D$, with the monotone convergence theorem, we deduce

$$\mathbb{E}\int_0^\tau \|\nabla v(X_t)\sigma(X_t)\mathbf{1}_{t<\tau}\|^2 \, dt < \infty.$$

Then

$$\mathbb{E}\int_0^\tau \|\nabla v(X_t)\sigma(X_t)\mathbf{1}_{t<\tau} - Z_t\|^2 \, dt = 0.$$

This achieves the proof of the proposition. □

From Lemma 6, we can deduce that $u_n$ belongs to $C^0(\overline{D}) \cap C^2(D)$ if $\sigma$ and $b$ belong to $C^1(D)$.

PROOF OF PROPOSITION 11. We fix $\eta > 0$ and we consider the set $D_\eta = D \setminus \{x \in D, \rho(x) \leq \eta\}$ for all $\eta > 0$. From Lemma 6, $u_n \in C^2(D)$ and satisfies $-\mathcal{L}u_n + u_n^{1+q} = 0$ on $D$. And on $D_\eta$, $u_n$ is bounded by $C/\eta^{2/q}$. Therefore, from Theorem 15.5 in [9], we obtain that $\nabla u_n$ is bounded on $D_{2\eta}$. The sequence is bounded in $C^1(D_{2\eta})$, thus the limit $u$ is continuous on $D_{2\eta}$, that is, $u$ is continuous on $D$. Moreover, we already know that $u$ is continuous on the boundary. Therefore, we deduce that $u$ belongs to $C^0(\overline{D}, \overline{\mathbb{R}}_+)$.

Now if we consider the PDE,

$$-\mathcal{L}v - v|v|^q = 0 \quad \text{on } D_\eta,$$
$$v = u \quad \text{on } \partial D_\eta,$$

from Theorem 15.18 in [9], the equation has a regular solution $v \in C^0(\overline{D_\eta}) \cap C^2(D_\eta)$. But this solution is also a continuous viscosity solution. Since $u$ is now a continuous viscosity solution of the same PDE, from the comparison result in [5], we deduce that $v = u$, that is, $u \in C^2(D_\eta)$. Hence, $u$ belongs to $C^2(D)$. □

Now we want to prove that $u$ is continuous on $D$ without the regularity conditions on $b$ and $\sigma$ of the Proposition 11. We just assume that (M), (L), (E) and (B) hold.

PROPOSITION 12. *The viscosity solution $u$ is continuous on $\overline{D}$ and is locally Hölder continuous on $D$.*

PROOF. We will show that for all open sets $D' \subset D$ such that $\overline{D'} \subset D$, there exists $0 < \alpha < 1$ such that the sequence of functions $u_n$ is bounded in the space $C^\alpha(D')$. $C^\alpha(D')$ is the set of functions $v$ such that

$$\|v\|_\alpha = \sup\left\{\frac{|v(x) - v(y)|}{|x - y|^\alpha},\ (x, y) \in D'\right\} < +\infty.$$

Since $u_n$ converges to $u$, we deduce that $u$ belongs to $C^\alpha(D')$ and thus is continuous on $D$.

In order to prove that $u_n$ is a bounded sequence in $C^\alpha(D')$, we will construct a sequence $v_m$ which will belong to $C^\alpha(D')$ and such that there exists a constant $K$ such that, for all $m \in \mathbb{N}$, $\|v_m\|_\alpha \leq K$. Let $b_m$ and $\sigma_m$ be two sequences of functions such that:

1. $b_m$ and $\sigma_m$ belong to $C^1(D)$ and $b_m$ and $\sigma_m$ are bounded in $L^\infty(D)$;
2. $b_m$ (resp. $\sigma_m$) converges to $b$ (resp. $\sigma$), uniformly on $D$;

3. $\sigma_m$ satisfies the condition (E).

Let $v_m$ be the unique solution in $C^0(\overline{D}) \cap C^2(D)$ (see Lemma 6 or [9]) of the equation

(6)
$$-\mathcal{L}_m v_m + v_m |v_m|^q = 0 \quad \text{on } D,$$
$$v_m = g \wedge n \quad \text{on } \partial D,$$

where $\mathcal{L}_m$ is the operator:

$$\forall x \in \mathbb{R}^d \quad \mathcal{L}_m \varphi(x) = \tfrac{1}{2} \operatorname{Trace}(\sigma_m \sigma_m^*(x) D^2 \varphi(x)) + b_m(x) \nabla \varphi(x).$$

For $x \in \overline{D}$, let $X^{x,m}$ be the solution of the SDE

$$\forall t \geq 0 \quad X_t^{x,m} = x + \int_0^t b_m(X_r^{x,m}) \, dr + \int_0^t \sigma_m(X_r^{x,m}) \, dB_r,$$

$\tau_m$ is the first exit time from $\overline{D}$ of the diffusion $X^{x,m}$, $(Y^{x,n,m}, Z^{x,n,m})$ is the solution of the BSDE:

$$Y_t^{x,n,m} = (g \wedge n)(X_{\tau_m}^{x,m}) + \int_t^{\tau_m} Y_r^{x,n,m} |Y_r^{x,n,m}|^q \, dr - \int_t^{\tau_m} Z^{x,n,m} \, dB_r.$$

From classical results on the SDE, $X^{x,m}$ converges to $X^x$ solution of the SDE (4) and the process $(Y^{x,n,m}, Z^{x,n,m})$ converges to $(Y^{x,n}, Z^{x,n})$ solution of the BSDE (33) (see Proposition 4.4 in [6]).

From Lemma 6, we have

$$v_m(x) = Y_0^{x,n,m} \quad \text{and} \quad u_n(x) = Y_0^{x,n}.$$

Therefore, $v_m$ converges to $u_n$. Moreover, we know that $v_m$ is a bounded sequence in $L^\infty(D)$.

Let $D'$ be a open subset of $D$ such that $\overline{D'} \subset D$. We apply Theorem 8.24 in [9]. The function $v_m$ is the solution of

$$\mathcal{L} v_m = v_m |v_m|^q = g \in L^\infty.$$

Therefore, there exists a real $0 < \alpha < 1$ and a constant $K$ such that

$$\|v_m\|_\alpha \leq K \|v_m\|_{L^\infty}.$$

The constants depend on the ellipticity constant of $\sigma_m$, on the bound on $b_m$ and $\sigma_m$ in $L^\infty$ and on the distance between $D'$ and $\partial D$. We deduce that $u_n$ belongs to $C^\alpha(D')$ and the norm $\|u_n\|_\alpha$ is bounded w.r.t. $n \in \mathbb{N}$.

Finally, $u$ belongs to $C^\alpha(D')$. □

5.3. *Minimal viscosity solution.* We prove the following:

THEOREM 9. *If $v$ is another nonnegative viscosity solution of the PDE (6) (in the sense of Definition 4), and if $v_* \geq g$ on $\partial D$, then $u \leq v$ on $\overline{D}$.*

PROOF. We show that for all $n \in \mathbb{N}^*$, $u_n \leq v_*$. The proof is the same as the proof of Theorem 3.3 in [5]. We fix $n$, we assume that there exists $z \in \overline{D}$ such that $\delta = u_n(z) - v_*(z) > 0$ and we will find a contradiction. The main tool is Theorem 3.2 in [5]. □

**6. Other generators $f$.** We have considered the generator $f(y) = -y|y|^q$. The main properties of this function are it is nonincreasing and allows the explosion at time $\tau$ [see (5) for the definition of this stopping time]. But we can also consider more general generators. Let $f : \mathbb{R} \to \mathbb{R}$ be a nonincreasing function of class $C^1$, such that there exists $q > 0$, $\kappa > 0$ s.t.

$$\forall y \geq 0 \qquad f(y) \leq -\kappa y^{1+q}. \tag{9}$$

The BSDE (10) has a unique solution if $\xi$ satisfies the condition (H6). From Remark 4, if $\xi \in L^\infty$, then (H6) holds. We also assume that $f(0) = 0$; thus, if $\xi \geq 0$, then $Y_t \geq 0$ for all $t \geq 0$.

First of all, the conclusion of Theorem 6 holds: there exists a constant $C$ such that for every solution $(Y, Z)$ of the BSDE (10),

$$\forall t \geq 0 \qquad |Y_t| \leq \frac{C}{(\rho(X_{t \wedge \tau}))^{2/q}}.$$

Indeed, with the notation of the proof of Theorem 6, equality (13) becomes

$$\Psi(X_{t \wedge \tau}) = \Psi(X_{T \wedge \tau}) + \int_{t \wedge \tau}^{T \wedge \tau} f(\Psi(X_r))\, dr$$
$$- \int_{t \wedge \tau}^{T \wedge \tau} \nabla \Psi(X_r) \sigma(X_r)\, dB_r$$
$$- \int_{t \wedge \tau}^{T \wedge \tau} [\nabla \Psi(X_r) b(X_r)$$
$$+ \tfrac{1}{2} \operatorname{Trace}(\sigma \sigma^*(X_r) D^2 \Psi(X_r)) + f(\Psi(X_r))]\, dr.$$

With assumption (9), equation (14) becomes

$$-(\nabla \Psi) b - \frac{1}{2} \operatorname{Trace}(\sigma \sigma^* D^2 \Psi) - f(\Psi)$$
$$\geq C \theta^{-2/q-2} \bigg[ \kappa C^q + \frac{2\theta}{q} (\nabla \theta) b$$
$$- \frac{1}{q}\bigg(\frac{2}{q} + 1\bigg) \|\sigma \nabla \theta\|^2 + \frac{\theta}{q} \operatorname{Trace}(\sigma \sigma^* D^2 \theta) \bigg];$$

and we can choose $C$ such that the right-hand side is nonnegative. The rest of the proof remains the same.

Now we suppose that $\xi$ is a nonnegative, $\mathcal{F}_\tau$-measurable random variable such that $\mathbf{P}(\xi = +\infty) > 0$. As in Section 2, we construct a process $(Y, Z)$ satisfying the conditions (D1) and (D2) of Definition 2: $(Y, Z)$ is the limit of the sequence of solutions $(Y^n, Z^n)$ of the BSDE (10) with terminal condition $\xi \wedge n$.

For the continuity of $Y$ [condition (D3) of Definition 2], Section 3.2 remains unchanged, if we have already proved that the limit of $Y_{t \wedge \tau}$ when $t$ goes to $+\infty$ exists.

We define the following function $F$ on $\mathbb{R}_+^*$:

$$F(y) = -\int_y^{+\infty} \frac{1}{f(x)} \, dx.$$

With (9) and since $f$ is of class $C^1$ and nonincreasing, $F$ is a positive, decreasing and convex function such that $\lim_{y \to 0} F(y) = +\infty$, and $\lim_{y \to +\infty} F(y) = 0$. Moreover, for $\varepsilon > 0$, $F_\varepsilon$ is defined by $F_\varepsilon(y) = F(y + \varepsilon)$.

Now if $\alpha$ is a constant such that $\xi \geq \alpha > 0$, for all $n \in \mathbb{N}$, all $t \geq 0$, $Y_t^n > 0$ a.s. We can apply the Itô formula to $F_\varepsilon(Y^n)$: for all $0 \leq t \leq T$,

$$(39) \quad \begin{aligned} F_\varepsilon(Y_{t \wedge \tau}^n) &= \mathbb{E}^{\mathcal{F}_t} F_\varepsilon(Y_{T \wedge \tau}^n) + \mathbb{E}^{\mathcal{F}_t} \int_{t \wedge \tau}^{T \wedge \tau} F_\varepsilon'(Y_r^n) f(Y_r^n) \, dr \\ &\quad - \tfrac{1}{2} \mathbb{E}^{\mathcal{F}_t} \int_{t \wedge \tau}^{T \wedge \tau} F_\varepsilon''(Y_r^n) \|Z_r^n\|^2 \, dr. \end{aligned}$$

Now $F_\varepsilon'' \geq 0$ and $0 \leq F_\varepsilon'(y) f(y) = \frac{f(y)}{f(y+\varepsilon)} \leq 1$. Let $T$ go to $+\infty$:

$$\begin{aligned} F_\varepsilon(Y_{t \wedge \tau}^n) &\leq \mathbb{E}^{\mathcal{F}_t} F_\varepsilon(\xi \wedge n) + \mathbb{E}^{\mathcal{F}_t}(\tau - t \wedge \tau) \leq F_\varepsilon(\alpha) + \mathbb{E}^{\mathcal{F}_t}(\tau - t \wedge \tau) \\ &\leq F(\alpha) + \mathbb{E}^{\mathcal{F}_t}(\tau - t \wedge \tau). \end{aligned}$$

Recall that from (C2), $\tau \in L^1(\Omega)$. Hence, for all $t \geq 0$, $\sup_n F(Y_t^n)$ belongs to $L^1$. With (39) and the same ideas as in the proof of Proposition 6, we deduce that, for $t \geq 0$,

$$F(Y_{t \wedge \tau}) = \mathbb{E}^{\mathcal{F}_t}(F(\xi) + \tau - t \wedge \tau) - \Phi_t,$$

where $\Phi$ is a nonnegative supermartingale.

Finally, if $f$ is a nonincreasing and $C^1$ function with $f(0) = 0$, such that (9) holds, and if $\xi$ is a nonnegative, $\mathcal{F}_\tau$-measurable random variable such that $\mathbf{P}(\xi = +\infty) > 0$, then the BSDE (10) has a minimal solution (in the sense of Definition 2). And the associated PDE (11) has also a minimal viscosity solution.

The only thing which we cannot describe just with inequality (9), is the behavior of $Y$ on the set $\{\xi = +\infty\}$ (see Proposition 1).

**Acknowledgments.** The author wishes to thank the referee and Professor Etienne Pardoux for the attention they paid to this article.

L.A.T.P. 39 RUE F. JOLIOT CURIE
13453 MARSEILLE CEDEX 13
FRANCE
E-MAIL: popier@cmap.polytechnique.fr
URL: http://www.cmap.polytechnique.fr/~popier